\documentclass[preprint,12pt,3p]{elsarticle}
\usepackage{amssymb}
\usepackage{color}
\usepackage{graphicx}
\usepackage{amsmath}
\usepackage{caption}
\usepackage{subcaption}
\usepackage{float}
\usepackage{hyperref}
\usepackage{textcomp}
\usepackage{amssymb}
\usepackage{enumitem} 
\newtheorem{theorem}{Theorem}[section]
\newtheorem{corollary}{Corollary}[section]
\newtheorem{definition}{Definition}[section]

\newtheorem{example}{Example}[section]
\newtheorem{note}{Note}[section]
\newtheorem{remark}{Remark}[section]

\journal{}

\begin{document}
\begin{frontmatter}

\title{Analytical aspects of matrix interpolation problems and its applications}

\author{Dharm Prakash Singh}%\fnref{label3}}
%\address{Department of Applied Mathematics, School of Vocational Studies and Applied Sciences, \\ Gautam Buddha University, Gr. Noida - 201312, India.}
%\address[label2]{Address Two\fnref{label4}}

%\author[label1,label2]{Author One\corref{cor1}\fnref{label3}}
%\address[label1]{Address One}
%\address[label2]{Address Two\fnref{label4}}

%cortext[cor1]{I am corresponding author}
%\fntext[label3]{I also want to inform about\ldots}
%\fntext[label4]{Small city}

\ead{prakashdharm@gmail.com}
%\ead[url]{author-one-homepage.com}

\author{Amit Ujlayan}
\address{Department of Applied Mathematics, School of Vocational Studies and Applied Sciences, \\ Gautam Buddha University, Gr. Noida - 201312, India.}
%\address[label5]{School of Vocational Studies and Applied Sciences  \\ Gautam Buddha University, Gr. Noida - 201310, India.}
\ead{iitk.amit@mail.com}

%\footnote{The work of first author is partially supported by the University Grant Commission, Ministry of Human Resource Development, New Delhi, India (JRF Award Letter No. F.16-6(DEC.2016)/2017(NET)).}
%\author[label5]{Author Two}
%\address[label5]{Some University}
%\ead{author.two@mail.com}

%\author[label1,label5]{Author Three}
%\ead{author.three@mail.com}

\begin{abstract}
In this paper, the $mn$-dimensional space of tensor-product polynomials of two variables, of degree at most $(m-1)+(n-1)$, is considered. A theory of two-variate polynomials is developed by establishing the algebra and basic algebraic properties with respect to the usual addition, scalar multiplication, and a newly defined algebraic operation in the considered space. Further, the existence of the considered space is established with respect to the matrix interpolation problem (MIP), $P(i,j)=a_{ij}$ for all $1 \leq i \leq m$, $1 \leq j \leq n$, corresponds to a given matrix $(a_{ij})_{m \times n}$ in the space of $m \times n$ order real matrices. The poisedness of the MIP is proved and three formulae are presented to construct the respective polynomial in the considered space. After that, using construction formulae, a polynomial map from the space of $m \times n$ order real matrices to the considered space is defined. Some properties of the polynomial map are investigated and some isomorphic structures between the spaces are installed. It is proved that the considered space is isomorphic to the space of $m \times n$ order real matrices with respect to the algebra structure. The polynomials in the considered space with respect to the MIP's for the given matrices also preserve the geometric properties of the matrices such as transpose, symmetry, and skew-symmetry. Some examples are included to demonstrate and verify the results. 
\end{abstract}

\begin{keyword}
Two-variate polynomial interpolation \sep polynomial algebra \sep algebraic structures 

{Mathematics Subject Classification (2010):} 41A05 \sep 41A10 \sep 41A63 \sep 65F99 \sep 15A99 
\end{keyword}
\end{frontmatter}

% \linenumbers

\section{Introduction}
Let $\mathbb{R}^{m \times n}$ denotes the space of $m \times n$ order real matrices and $(a_{ij})_{m \times n}$ denotes the matrix $A=(a_{ij})_{m \times n}$, where $m,n \in \mathbb{N}$, \cite{Alfio,axler15}. Again, let $\Pi^{d}$ be the space of all $d$-variate polynomials and $\Pi_{N}^{d}$ be  ${{d+N}\choose{N}}$-dimensional subspace of $\Pi^{d}$, of degree at most $N$, over the field $\mathbb{F}$, where $d,N \in \mathbb{N}$, \cite{carl:1,Yuan:1}. For a given finite linearly independent set $\Theta$ of functionals and an associated vector $Y = (y_{\theta} : \theta \in \Theta)$ of prescribed values, find a polynomial $p$ in $\Pi^{d}$ such that
\begin{equation}\label{general interpolation problem}
    \Theta p = Y,   \quad   i.e., \quad     \theta p = y_{\theta}, \quad \theta \in \Theta,
\end{equation}
represents general form of the polynomial interpolation problem, \cite{Thomas:1,carl:1}. 

Let $\mathcal{S}$ be a subspace of $\Pi^{d}$, the polynomial interpolation problem \eqref{general interpolation problem}
with respect to the set $\Theta$ is said to be \emph{poised} or \emph{correct} for $\mathcal{S}$, if for any $Y \in \mathbb{F}^{\Theta}$ there exists a unique $p \in \mathcal{S}$ such that $\Theta p = Y$, \cite{Thomas:1}. In other words, the interpolation problem \eqref{general interpolation problem} with respect to a finite set $\Theta$ of the cardinality $q$, $q \in \mathbb{N}$ is \textit{poised} in $\mathcal{S}$, if and only if, $q=dim(\mathcal{S})$ and the determinant of the sample matrix, \cite{Kamron:1,Olver:1},
\begin{equation}\label{sample matrix}
\Theta \mathcal{B}=\begin{bmatrix} \theta p: \theta \in \Theta, p \in \mathcal{B} \end{bmatrix} \in \mathbb{F}^{{\Theta}\times \mathcal{B}}
\end{equation}
is non-zero for any choice of basis $\mathcal{B}$ of $\mathcal{S}$. The problem is said to be \emph{singular} if $det(\Theta \mathcal{B})=0$ for every choice of the set $\Theta$, \emph{regular} if $det(\Theta \mathcal{B})=0$ for no choice of the set $\Theta$, and almost regular if $det(\Theta \mathcal{B})=0$ only on a subset of $\mathbb{R}^{d}$ of measure zero, \cite{Thomas:1,Kamron:1,Olver:1,Gasca:1,Mehaute:1}. 

%For the given set $\Theta$ and a subspace $\mathcal{H} \subset \Pi^d$, if the problem \eqref{lag:int} is poised in $\mathcal{H}$, then the construction of the polynomial $p \in \mathcal{H}$ which satisfy \eqref{lag:int}, can be completed using the Lagrange or Newton  interpolation formulae discussed in \cite{Yuan:1,Kamron:1}.

\subsection{Motivation and the considered problem}
Nowadays, a matrix is the majorly used data structure in computer-oriented models and algorithms of science and engineering. Particularly, in computer graphics \cite{image:2}, signal processing \cite{image:1}, image processing \cite{image:2,image:1}, computer-aided designs in control systems \cite{CAD:1,CAD:2}, thermal engineering \cite{thermal:1} to name a few. The generation of the respective three-dimensional smooth interpolating surface is one of the important tools for analysis, comparison, and interpretation. In medical imaging, the surface interpolation plays a very important role in computer-assisted surgical planning and medical diagnosis, \cite{surface:1}.

The matrices and matrix algebra are used in a very large-scale to solve the computational problems as their solutions constitute an important and major section of scientific, engineering, and numerical computations on the modern computers. The computational problems build a challenging and application-oriented segment of effectual research. In last decades, an intense research has been witnessed to develop efficient numerical algorithms for the matrix computations and a rapid progress have been noticed in the formation of efficient and effective algorithms for the polynomials in support of algebraic, symbolic, numerical, and scientific computations, \cite{imp:2,imp:1,imp:3}.

In this regard, considering the matrix interpolation problem \cite{dpmip:1}, defined as follows:
\begin{definition} Let $(a_{ij}) \in \mathbb{R}^{m \times n}$ be a given matrix. For a given subspace $\mathcal{P}$ of $\Pi^{2}$ and the set $\mathcal{D}=\left\{(i,j):i=1,2,\ldots,m, j=1,2,\ldots,n \right \}$ of pairwise distinct interpolation points or nodes, find a polynomial $P \in \mathcal{P}$, such that 
\begin{equation}\tag {*}
            P(i,j)=a_{ij} \textrm{ for all }i=1,2,\ldots,m,\;j=1,2,\ldots,n \label{mip}
	\end{equation}
\noindent is defined as \emph{Matrix Interpolation Problem  (MIP)}.
\end{definition}

For the given matrix $(a_{ij}) \in \mathbb{R}^{m \times n}$, the MIP \eqref{mip} \textit{can not be poised} in $\Pi_{N}^{2}$ for almost every $N \in \mathbb{N}$, in view of the fact that,
\begin{enumerate}[label=(\roman*)]
\item if $mn=dim\Pi_{N}^{2}$, the determinant of the respective sample matrix is zero.
\item if $mn \neq dim\Pi_{N}^{2}$, the necessary condition of $\frac{(N+1)(N+2)}{2}$ interpolation points does not hold (see section 1.1 in \cite{dpmip:1}).
\end{enumerate}

%\begin{remark} The construction of the sample matrix of the order $\frac{(N+1)(N+2)}{2} \times \frac{(N+1)(N+2)}{2}$ corresponds to the MIP \eqref{mip} such that $mn = dim\Pi_{N}^2$ is analytically intractable for large $N$. However, using the suitable computer program, it is verified that the determinant of such sample matrices is zero for $N \leq 30$ (higher dimension required more execution time). We conjecture that it is true for all $N \in \mathbb{N}$.
%\end{remark}

\subsection{Brief review of literature}
Let $k$ be a positive integer and $ \Theta = \left\{x_1, x_2, \ldots, x_k \right\} \subset \mathbb{R}^d$ be the set of pairwise disjoint points. Then, for a given subspace $\mathcal{H} \subset \Pi^d$ and some given constants $f_1, f_2, \ldots, f_k$ in $\mathbb{R}$, the problem to find a polynomial $p \in \mathcal{H}$ with respect to $\Theta$, such that 
	\begin{eqnarray}\label{lag:int}
            p(y_i)=f_i \textrm{ for all } i=1,2,...,k 
	\end{eqnarray}
can be defined as Lagrange interpolation problem. The interpolation points $x_i$, $i=1,2,...,k$ are called interpolation sites or nodes and $\mathcal{H}$ is known as the interpolation space, \cite{carl:1, Gasca:1}. 

For the given set $\Theta$ and a subspace $\mathcal{H}$, if the problem \eqref{lag:int} is poised in $\mathcal{H}$, then the required polynomial $p \in \mathcal{H}$ can be constructed using the interpolation formulae discussed in \cite{Yuan:1,Kamron:1}.

In general, the problem \eqref{lag:int} need not be \emph{poised} in a given subspace $\mathcal{H}$ of $\Pi_{N}^{d}$, with respect to the given set $\Theta$ of dimension ${{d+N}\choose{N}}$. As in the case of $d \geq 2$, for a given ${{d+N}\choose{N}}$-dimensional subspace, of a certain total degree, the poisedness of the problem does not independent of the position and geometry (or configuration) of the data points. Also, the basic structure and nature of the interpolation changes utterly and the existence of \emph{Haar Spaces} of dimension higher than one vanishes, \cite{carl:1,Yuan:1,Thomas:1,Gasca:1}.

However, there always exist at least one set $\Theta$ of dimension ${{d+N} \choose{N}}$ in $\mathbb{R}^{d}$ such that the interpolation problem \eqref{lag:int} can poised in $\Pi_{N}^{d}$, \cite{Yuan:1,Thomas:1,Gasca:1}. Again, the Kergin interpolation \cite{kergin:1} insures that, there always exist at least one subspace of $\Pi_{r-1}^{d}$ in which the interpolation problem \eqref{lag:int} can be poised for some arbitrarily known set $\Theta$ of dimension $r$, $r \in \mathbb{N}$. Hence, the Lagrange interpolation problem \eqref{lag:int} can never be singular, \cite{Gasca:1}. The work in this manuscript focuses on a special case of two-variate interpolation problem. Therefore, we restrict our self to $d=2$.

For $d=2$, the MIP \eqref{mip} is a class of the Lagrange interpolation problems. Therefore, due to Kergin, there must exist at at least one $mn$-dimensional subspace of $\Pi^{2}$ in which the MIP \eqref{mip} can always be poised with respect to the set $\mathcal{D}$ for a given matrix $a_{ij} \in \mathbb{R}^{m \times n}$. In \cite{dpmip:1}, the constructive existence of two $mn$-dimensional subspaces $_{n}^{1}\Pi_{m,n}^2$ and $_{1}^{m}\Pi_{m,n}^2$ has been established, by projecting two-dimensional interpolating points into one-dimension such that all projection are different using a bijective transformation approach and univariate polynomial interpolation, in which the MIP \eqref{mip} can always be poised, where $_{\alpha}^{\beta}\Pi_{m,n}^{2}$ is the the space of all two variable polynomials in $x$ and $y$, of degree at most $mn-1$, in two parameters $\alpha$ and $\beta$, with real coefficients, such that, if $P\in {_{\alpha}^{\beta}\Pi_{m,n}^{2}}$, then
\begin{eqnarray}\label{dpsofk}
        P(x,y)= \sum_{k=0}^{mn-1}{{\lambda_{k}}(\alpha x + \beta y)^{k}}.
\end{eqnarray}

Further, it is concluded that the MIP (\ref{mip}) can always be \textit{poised} in $_{\alpha}^{\beta}\Pi_{m,n}^{2}$ for unaccountably infinite number of choices of the parameters $\alpha$ and $\beta$, provided the condition 
$\alpha i_1 + \beta j_1 \neq \alpha i_2 + \beta j_2$ for all $i_1 \neq i_2$ and $j_1 \neq j_2$ holds, where $i_1,i_2 \in \left \{1,2,\ldots,m \right \} \textrm{ and } j_1,j_2 \in \left \{1,2,\ldots,n\right \}$.

From the computational and algebraic point of view, we have identified some flaws of the solutions and the established subspaces which are listed as follows:
\begin{enumerate}
    \item The degree of the interpolating polynomials is up to $mn-1$, which is somehow high from the computational point of view.
    \item The interpolating polynomials in these subspaces, corresponds to the MIP \eqref{mip}, do not satisfy the geometric structures of the given matrix $a_{ij} \in \mathbb{R}^{m \times n}$ such as transpose, symmetry, and skew-symmetry.
    \item The suggested subspaces are isomorphic to $\mathbb{R}^{m \times n}$ with respect to the vector algebra structure, but they do not preserve multiplication structure as defined in the space of matrices $\mathbb{R}^{m \times n}$.
\end{enumerate}

The identified demerits in the solutions and the interpolation spaces is motivating us to identify at least one $mn$-dimensional space of two variable polynomials, of least possible degree, in which the MIP\eqref{mip} can always be poised and which can preserve almost every algebraic as well as geometric structures of the space $\mathbb{R}^{m \times n}$. 

In \cite{narumi1920some}, the author, Seimatsu Narumi suggested a $(m+1)(n+1)$-dimensional space of two variable polynomials, of total degree $mn$ (degree of the variables not greater than $m$ and $n$ respectively). The suggested space is known as the space of tensor-product polynomials in two variable. Narumi claimed this space as a correct interpolation space for a two-variate interpolation problem $P(a_i,b_j)=f(a_i,b_j)$, corresponds to a two-variate real-valued function $f(x,y)$, by sampling the data with respect to equally spaced interpolation points of the form $(a_i,b_j)$, $0 \leq i \leq m$, $0 \leq j \leq n$. A method of bivariate divided difference is discussed to evaluate the coefficients of the required polynomial in the suggested space, but the existence of the space or the approach of used interpolation method is not mentioned. 

For a given matrix $a_{ij} \in \mathbb{R}^{m \times n}$, the MIP \eqref{mip} becomes a special case of the Narumi interpolation problem, where $f:\mathbb{R}^2 \rightarrow \mathbb{R}$ and the interpolation points takes the form $(i,j)$, $1 \leq i \leq m$, $1 \leq j \leq n$. Therefore, in this part of work, the $mn$-dimensional space of tensor-product polynomials in two variable, of degree at most $(m-1)+(n-1)$, is considered.

\subsection{Objectives and arrangement of the article}
%The matrices are used in a very large-scale to solve the problems in scientific and engineering computing and a rapid progress has been noticed to improve the efficiency of the algorithms for the polynomial computations, in last decades. Therefore, 
The main objective of this paper is to develop a theory of two-variate polynomials in the considered space with respect to the algebraic operations and geometric properties of the space $\mathbb{R}^{m \times n}$. The second objective of this paper is to establish the existence of the considered space in which the MIP\eqref{mip} can be poised for all $a_{ij} \in \mathbb{R}^{m \times n}$ and present some formulae to construct the respective polynomial in the established space. Defining a polynomial map from $\mathbb{R}^{m \times n}$ to the considered space, investigation of some algebraic properties of the polynomial map, and establish some isomorphic structures between the spaces are the other prime objectives.

The rest of the paper is consolidated as follows. In section 2, the $mn$-dimensional space of tensor-product polynomials in two variables, of degree at most $(m-1)+(n-1)$, will be defined. The basic algebra and some algebraic properties of the defined space will be discussed with respect to usual addition, scalar multiplication, and a newly defined algebraic operation. In section 3, the existence of the defined space will be established and the poisedness of the MIP \eqref{mip} in it will be proved. Some formulae will be presented to construct the two-variate polynomial in the space corresponds to the MIP \eqref{mip} for a given matrix $(a_{ij})\in \mathbb{R}^{m \times n}$. In section 4, the polynomial map from $\mathbb{R}^{m \times n}$ to the established subspace will be defined. Some properties of the additive and multiplicative structures will be discussed. In section 5, some examples will be included to validate the theoretical findings and to observe the geometrical prospects. Finally, section 6 will conclude the paper.

%In \cite{Gasca:1}, \emph{Mariano Gasca} and \emph{Thomas Sauer} concluded that ``To date, very little can be said on the usefulness of polynomial interpolation in serious applications''. 

\section{A theory of two-variate polynomials}
In this section, the $mn$-dimensional space of tensor-product polynomial in two variables, of degree at most $(m-1)+(n-1)$, will be defined. Some properties of the defined space will be investigated and a theory of two-variate polynomials will be developed in the defined space with respect to some algebraic operations and geometric properties of the space $\mathbb{R}^{m \times n}$. Some definitions, examples, counterexamples, remarks, and theorems will be included to identify the properties of the defined space.

\begin{definition}
Let $\mathcal{P}_{m}^{n}$ denotes the space of two variable polynomials in $x$ and $y$, of degree at most $(m-1)+(n-1)$, with real coefficients such that, if $P\in \mathcal{P}_{m}^{n}$, then
    \begin{equation}\tag{**}
        P(x,y)= {\sum_{k_1=0}^{m-1} {\sum_{k_2=0}^{n-1}}}{\lambda _{{k_1,}{k_2}}}{x^{k_1}}{y^{k_2}}. \label{dpsosk} 
    \end{equation} 
Here, \textrm{dim }$\mathcal{P}_{m}^{n}=mn$.
\end{definition}

\subsection{Equality, addition, and scalar multiplication of the polynomials in $\mathcal{P}_{m}^{n}$}
The equality, addition and scalar multiplication of the polynomials in $\mathcal{P}_{m}^{n}$ is defined as the standard definitions of the two-variate polynomials.

\begin{theorem}\label{thm:vector space}
       The space $\mathcal{P}_{m}^{n}$ is an $mn$-dimensional subspace of $\Pi_{m+n-2}^{2}$.
\end{theorem}
\noindent \textbf{Proof:} The space $\mathcal{P}_{m}^{n}$ is an $mn$-dimensional non-empty subset of the vector space $\Pi_{m+n-2}^{2}$. Suppose $P, Q \in \mathcal{P}_{m}^{n}$ be two polynomials, given as
\begin{eqnarray}
       P(x,y)= {\sum_{k_1=0}^{m-1} {\sum_{k_2=0}^{n-1}}}{\lambda _{{k_1,}{k_2}}}{x^{k_1}}{y^{k_2}}, \quad Q(x,y)= {\sum_{k_1=0}^{m-1} {\sum_{k_2=0}^{n-1}}}{\mu _{{k_1,}{k_2}}}{x^{k_1}}{y^{k_2}}, \nonumber
\end{eqnarray}
and $\alpha$ be a real scalar. Then,
\begin{flalign*}
\qquad P(x,y)+Q(x,y) & = {\sum_{k_1=0}^{m-1} {\sum_{k_2=0}^{n-1}}}{\lambda _{{k_1,}{k_2}}}{x^{k_1}}{y^{k_2}} + {\sum_{k_1=0}^{m-1} {\sum_{k_2=0}^{n-1}}}{\mu _{{k_1,}{k_2}}}{x^{k_1}}{y^{k_2}} \\
        & = {\sum_{k_1=0}^{m-1} {\sum_{k_2=0}^{n-1}}}({\lambda_{{k_1,}{k_2}}+ \mu_{{k_1,}{k_2}}}){x^{k_1}}{y^{k_2}}, && 
\end{flalign*}
\begin{flalign*}
\textrm{and }\;\; \alpha P(x,y) = \alpha {\sum_{k_1=0}^{m-1}{\sum_{k_2=0}^{n-1}}}(\lambda_{{k_1,}{k_2}}){x^{k_1}}{y^{k_2}} = {\sum_{k_1=0}^{m-1} {\sum_{k_2=0}^{n-1}}}({\alpha}\lambda_{{k_1,}{k_2}}){x^{k_1}}{y^{k_2}}.&&
\end{flalign*}
i.e., $P(x,y)+Q(x,y), \alpha P(x,y) \in \mathcal{P}_{m}^{n}$ for some scalars $\delta_{{k_1,}{k_2}}={\lambda_{{k_1,}{k_2}}+ \mu_{{k_1,}{k_2}}}$ and $\sigma_{{k_1,}{k_2}}={\alpha}\lambda_{{k_1,}{k_2}}$. This completes the proof.

\begin{corollary}\label{coro:isomorphism1}
The spaces $\mathbb{R}^{m \times n}$ and $\mathcal{P}_{m}^{n}$ are isomorphic with respect to the vector space structure.
\end{corollary}

\begin{remark}
The set of standard basis of the vector space $\mathcal{P}_{m}^{n}$, is  
\begin{equation}
\left\{x^{k_1-1}y^{k_2-1}:1 \leq k_1 \leq m, 1 \leq k_2 \leq n, \textrm{ where } k_1,k_2,m,n \in \mathbb{N} \right\}.
\end{equation}
\end{remark}

\subsection{Product of the polynomials in $\mathcal{P}_{m}^{n}$}
\begin{definition}\label{definition:dp product}
Let $P \in \mathcal{P}_{m}^{n}$ and $Q \in \mathcal{P}_{n}^{q}$ be two polynomials, then the operation `$\otimes$' given by
    \begin{equation}\label{definition:multiplication}
        P(x,y)\otimes Q(x,y)= {\sum_{k=1}^{n}(P(x,k) \times Q(k,y))} 
        \textrm{ in } \mathcal{P}_{m}^{q},
    \end{equation}
is defined as DP product.
\end{definition}

\begin{example}\label{ex:p1}
Let $P \in \mathcal{P}_{2}^{2}$ and $Q \in \mathcal{P}_{2}^{3}$ be two polynomials given as $P(x,y)=14x-10xy+13y-18$ and $Q(x,y)=\frac{1}{2}(14y^2-11xy^2+41xy-34x-52y+44)$
respectively. Then, the definition \eqref{definition:dp product} implies that
    $P(x,y)\otimes Q(x,y) = \frac{1}{2}(60xy^2-224xy-79y^2+184x+295y-242)$ in  $\mathcal{P}_{2}^{3}$. 
But, $Q(x,y)\otimes P(x,y)$ is not defined.
\end{example}

\begin{remark}
    If $P(x,y)\otimes Q(x,y)$ is defined, it is not necessary that $Q(x,y)\otimes P(x,y)$ is also defined.
\end{remark}

\begin{example}\label{ex:p2}
Let $P \in \mathcal{P}_{1}^{2}$ and $Q \in \mathcal{P}_{2}^{1}$ be two polynomials given as $P(x,y)=2y-3$ and $Q(x,y)=-2x+3$ respectively. Then, the definition \eqref{definition:dp product} implies that $P(x,y)\otimes Q(x,y) = -2 \in \mathcal{P}_{1}^{1}$  and $Q(x,y)\otimes P(x,y) = -4xy+6x+6y-9 \in \mathcal{P}_{2}^{2}$.
%where $P(x,y)\otimes Q(x,y) \in \P_{1}^{1}$ and $Q(x,y)\otimes P(x,y) \in \P_{2}^{2}$.
\end{example}

\begin{remark}
If $P \in \mathcal{P}_{m}^{n}$ and $Q \in \mathcal{P}_{p}^{q}$, then $P(x,y)\otimes Q(x,y)$ and $Q(x,y)\otimes P(x,y)$ both are defined, if and only if $n=p$ and $m=q$.
\end{remark}

\begin{example}\label{ex:p3}
Let $P \in \mathcal{P}_{2}^{2}$ and $Q \in \mathcal{P}_{2}^{2}$ be two polynomials given as $P(x,y)=-x-y+3$ and $Q(x,y)=-2xy+3i+3j-4$ respectively. Then, the definition \eqref{definition:dp product} implies that
$P(x,y)\otimes Q(x,y) = y-x \in \mathcal{P}_{2}^{2}$ and $Q(x,y)\otimes P(x,y) = x-y \in \mathcal{P}_{2}^{2}$,
%where $P(x,y)\otimes Q(x,y)$, $Q(x,y)\otimes P(x,y)\in \mathcal{P}_{2}^{2}$. 
Here, $P(x,y)\otimes Q(x,y) \neq Q(x,y)\otimes P(x,y)$.
\end{example}

\begin{remark}
Let $P, Q \in \mathcal{P}_{n}^{n}$, then $P(x,y)\otimes Q(x,y)$ and $Q(x,y)\otimes P(x,y)$ both are defined, but it is not necessary that $P(x,y)\otimes Q(x,y)=Q(x,y)\otimes P(x,y)$.
\end{remark}

\begin{example}\label{ex:p4}
Let $P \in \mathcal{P}_{2}^{2}$ and $Q \in \mathcal{P}_{2}^{2}$ be two polynomials given as $P(x,y)=3xy-3x-4y+4$ and $Q(x,y)=2xy-5x-4y+10$ respectively. Then, the definition \eqref{definition:dp product} implies that
$P(x,y)\otimes Q(x,y) =  0 \in \mathcal{P}_{2}^{2}$ and $Q(x,y)\otimes P(x,y) = xy-x-2y \in \mathcal{P}_{2}^{2}$. 
\end{example}

\begin{remark}
$P(x,y)\otimes Q(x,y)=0$ need not implies that either $P(x,y)=0$ or $Q(x,y)=0$ or $Q(x,y)\otimes P(x,y)=0$, i.e., the set $\mathcal{P}_{n}^{n}$ possess zero-divisor property.
\end{remark}

\subsubsection{Some properties of the DP product}
\begin{theorem}[DP product closure]\label{thm:multiplicative closure}
       The set $\mathcal{P}_{m}^{n}$ satisfies closure property under the operation $\otimes$, if and only if $m=n$.
\end{theorem}
\noindent \textbf{Proof:} Suppose $\mathcal{P}_{m}^{n}$ satisfies closure property under the operation $`\otimes '$, then it has to prove that $m=n$ for all $m,n \in \mathbb{N}$. Let $P_1, P_2 \in \mathcal{P}_{m}^{n}$ and if $m \neq n$, then by the definition of DP product \eqref{definition:dp product}, neither $P_1(x,y)\otimes P_2(x,y)$ nor $P_2(x,y)\otimes P_1(x,y)$ is defined. The contra-positive approach of proof, completes the proof of first part.

Conversely, suppose that $m=n$ for all $m,n \in \mathbb{N}$. Let us consider that $P_1, P_2 \in \mathcal{P}_{n}^{n}$, then $P_1(x,k) \in \mathcal{P}_{n}^{1}$ and $P_2(k,y) \in \mathcal{P}_{1}^{n}$, for some fixed $k$, where $1 \leq k \leq n$, $n \in \mathbb{N}$. Thus,
\begin{flalign*}
& \quad P_1(x,k) \times P_2(k,y) \in \mathcal{P}_{n}^{n} \textrm{ for all } 1 \leq k \leq n  \quad   \quad  \textrm{by definition of DP product}\\
& \Rightarrow {\sum_{k=1}^{n}(P_1(x,k) \times P_2(k,y))} \in \mathcal{P}_{n}^{n}  \quad   \quad   \quad   \quad  \quad   \quad   \quad \textrm{by additive closure}\\
& \Rightarrow P_1(x,y)\otimes P_2(x,y)) \in \mathcal{P}_{n}^{n}  \quad   \quad   \quad   \quad  \quad   \quad   \quad \quad \quad \textrm{by definition of DP product}&&
\end{flalign*}
This completes the proof.

\begin{theorem}[DP product and zero polynomial]\label{thm:zero polynomial}
       Suppose $P \in \mathcal{P}_{m}^{n}$, then
\begin{enumerate}[label=(\roman*)]
    \item $P(x,y)\otimes O(x,y)=0 \in \mathcal{P}_{m}^{p}$ for all $O(x,y)=0 \in \mathcal{P}_{n}^{p}$,
    \item $O(x,y)\otimes P(x,y)=0 \in \mathcal{P}_{p}^{n}$ for all $O(x,y)=0 \in \mathcal{P}_{p}^{m}$.
\end{enumerate}
\end{theorem}

\noindent \textbf{Proof:} Suppose $P \in \mathcal{P}_{m}^{n}$, then\\
(i) For all $O(x,y)=0 \in \mathcal{P}_{n}^{p}$, we get
\begin{equation}
P(x,y)\otimes O(x,y) = {\sum_{k=1}^{n}(P(x,k) \times O(k,y))} = {\sum_{k=1}^{n}(P(x,k) \times 0} = 0 \in \mathcal{P}_{m}^{p}. \nonumber 
\end{equation}
\noindent (ii) For all $O(x,y)=0 \in \mathcal{P}_{p}^{m}$, we get
\begin{equation}
O(x,y)\otimes P(x,y) = {} {\sum_{k=1}^{n}(O(x,k) \times P(k,y))} = {\sum_{k=1}^{n}(0 \times P(x,k))} = 0 \in \mathcal{P}_{p}^{n}. \nonumber 
\end{equation}
So every entry of the product is the scalar zero, i.e., the result is a zero polynomial. This completes the proof.

\begin{theorem}[DP product and scalar multiplication]\label{thm: scalar and polynomial multiplication}
Suppose $P \in \mathcal{P}_{m}^{n}$, $Q \in \mathcal{P}_{n}^{p}$, and $\alpha$ be a real scalar. Then, $$\alpha (P(x,y)\otimes Q(x,y))={\left(\alpha P(x,y)\right)}\otimes Q(x,y)=P(x,y)\otimes {\left(\alpha Q(x,y)\right)} \textrm{ in } \mathcal{P}_{m}^{p}.$$
\end{theorem}
\noindent \textbf{Proof:} Let $P \in \mathcal{P}_{m}^{n}$, $Q \in \mathcal{P}_{n}^{p}$ and $\alpha$ be a real scalar. Then,
\begin{flalign*}
\alpha (P(x,y)\otimes Q(x,y)) & = \alpha \left({\sum_{k=1}^{n}{\left(P(x,k) \times Q(k,y)\right)}}\right) \quad  \textrm{by definition of DP product}\\
    & =\alpha {\sum_{k=1}^{n}{\left(P(x,k) \times Q(k,y)\right)}} \quad  \textrm{by scalar multiplication}\\
    & = {\sum_{k=1}^{n}{\alpha \left(P(x,k) \times Q(k,y)\right)}} \quad  \textrm{by scalars distribution}\\
    & = {\sum_{k=1}^{n}{\left(\left(\alpha P(x,k)\right) \times Q(k,y)\right)}} \quad  \textrm{by scalars multiplication}\\
    & = {\left(\alpha P(x,y)\right)}\otimes Q(x,y) \quad  \textrm{by definition of DP product}&& 
\end{flalign*}
so the polynomials $\alpha (P(x,y)\otimes Q(x,y))$ and ${\left(\alpha P(x,y)\right)}\otimes Q(x,y)$ are equal, entry-by-entry, and by the definition of polynomial equality, they are equal polynomials in $\mathcal{P}_{m}^{p}$. Similarly, the polynomials $\alpha (P(x,y)\otimes Q(x,y))$ and $P(x,y) \otimes (\alpha Q(x,y))$ are equal in $\mathcal{P}_{m}^{p}$. This completes the proof.

\begin{theorem}[DP product associativity]\label{thm:multiplicative associativity}
Suppose $P \in \mathcal{P}_{m}^{n}$, $Q \in \mathcal{P}_{n}^{p}$ and $R \in \mathcal{P}_{p}^{q}$, then $$P(x,y)\otimes \left(Q(x,y)\otimes R(x,y)\right)=\left(P(x,y)\otimes Q(x,y)\right)\otimes R(x,y) \textrm{ in } \mathcal{P}_{m}^{q}.$$
\end{theorem}
\noindent \textbf{Proof:} Let $P \in \mathcal{P}_{m}^{n}$, $Q \in \mathcal{P}_{n}^{p}$ and $R \in \mathcal{P}_{p}^{q}$. Then,
%\begin{flalign*}
$P(x,y)\otimes \left(Q(x,y)\otimes R(x,y)\right)$
%\end{flalign*}
\begin{flalign*}
&= P(x,y) \otimes \left({\sum_{l=1}^{p}{\left(Q(x,l) \times R(l,y)\right)}}\right) \quad  \textrm{by definition of DP product}&&
\end{flalign*}
\begin{flalign*}
&= {\sum_{k=1}^{n}{\left(P(x,k) \times \left({\sum_{l=1}^{p}{\left(Q(k,l) \times R(l,y)\right)}}\right)\right)}} \quad  \textrm{by definition of DP product}&&
\end{flalign*}
\begin{flalign*}
& ={\sum_{k=1}^{n}{\sum_{l=1}^{p}}{P(x,k) \times Q(k,l) \times R(l,y)}} \quad  \textrm{by distributivity in $\mathbb{R}$}&&
\end{flalign*}
\begin{flalign*}
& ={\sum_{l=1}^{p}{\sum_{k=1}^{n}}{P(x,k) \times Q(k,l) \times R(l,y)}} \quad  \textrm{by switching the order of finite sums}&&
\end{flalign*}
\begin{flalign*}
& ={\sum_{l=1}^{p}{\sum_{k=1}^{n}}{R(l,y)\times P(x,k)\times Q(k,l)}} \quad  \textrm{by commutativity in $\mathbb{R}$}&&
\end{flalign*}
\begin{flalign*}
&= {\sum_{l=1}^{p}{\left(R(l,y) \times \left({\sum_{k=1}^{n}{\left(P(x,k) \times Q(k,l)\right)}}\right)\right)}} \quad \textrm{since $R(l,y)$ is independent of the index $k$} &&
\end{flalign*}
%\begin{flalign*}
%& \quad\quad\quad\quad\quad\quad\quad\quad\quad\quad  \textrm{since $R(l,y)$ is independent of the index $k$}&&
%\end{flalign*}
\begin{flalign*}
&= {\sum_{l=1}^{p}{\left(R(l,y) \times \left(P(x,l)\otimes Q(x,l)\right)\right)}} \quad  \textrm{by definition of DP product}&&
\end{flalign*}
\begin{flalign*}
&= {\sum_{l=1}^{p}{\left(\left(P(x,l)\otimes Q(x,l)\right) \times R(l,y) \right)}} \quad  \textrm{by commutativity in $\mathbb{R}$}&&
\end{flalign*}
\begin{flalign*}
&=\left(P(x,y)\otimes Q(x,y)\right)\otimes R(x,y) \quad  \textrm{by definition of DP product}&&
\end{flalign*}
So the polynomials $P(x,y)\otimes \left(Q(x,y)\otimes R(x,y)\right)$ and $\left(P(x,y)\otimes Q(x,y)\right)$ $\otimes R(x,y)$ are equal, entry-by-entry, and by the definition of polynomial equality, they are equal polynomial in $\mathcal{P}_{m}^{q}$. This completes the proof.

\begin{theorem}[DP product distributivity across addition]\label{thm:distributivity across addition}
Suppose $P,S \in \mathcal{P}_{m}^{n}$ and $Q,R \in \mathcal{P}_{n}^{p}$, then \\
(i) \textrm{Left Distributive Law:} $$P(x,y)\otimes {\left(Q(x,y)+R(x,y) \right)}=P(x,y)\otimes Q(x,y)+P(x,y)\otimes R(x,y) \textrm{ in } \mathcal{P}_{m}^{p}.$$ 
(ii) \textrm{Right Distributive Law:} $${\left(P(x,y)+S(x,y) \right)}\otimes Q(x,y)=P(x,y)\otimes Q(x,y)+S(x,y)\otimes Q(x,y) \textrm{ in } \mathcal{P}_{m}^{p}.$$
\end{theorem}
\noindent \textbf{Proof:} Let us consider that $P \in \mathcal{P}_{m}^{n}$, $Q \in \mathcal{P}_{n}^{p}$, $R \in \mathcal{P}_{n}^{p}$ and $S \in \mathcal{P}_{m}^{n}$. Then, 
\begin{flalign*}
(i) \textrm{ } P(x,y)\otimes {\left(Q(x,y)+R(x,y) \right)}
%\end{flalign*}
%\begin{flalign*}
&= {} {\sum_{k=1}^{n}(P(x,k) \times {\left(Q(k,y)+R(k,y) \right)}})\\
& = {} {\sum_{k=1}^{n}{\left(P(x,k) \times Q(k,y)+P(x,k) \times R(k,y) \right)}}\\
& = {} {\sum_{k=1}^{n}{\left(P(x,k) \times Q(k,y)\right)}} + {\sum_{k=1}^{n}{\left(P(x,k) \times R(k,y)\right)}}\\
& = {} P(x,y)\otimes Q(x,y)+P(x,y)\otimes R(x,y). && 
\end{flalign*}
So the polynomials $P(x,y)\otimes {\left(Q(x,y)+R(x,y) \right)}$ and $P(x,y)\otimes Q(x,y)$ + $P(x,y)$ $\otimes R(x,y)$ are equal, entry-by-entry, and by the definition of polynomial equality, they are equal polynomials in $\mathcal{P}_{m}^{p}$.
\begin{flalign*}
(ii) \textrm{ } {\left(P(x,y)+S(x,y) \right)}\otimes Q(x,y)
& = {} {\sum_{k=1}^{n}( {\left(P(x,k)+S(x,k)\right) \times Q(k,y))}}\\
    & = {} {\sum_{k=1}^{n}{\left(P(x,k) \times Q(k,y)+S(x,k) \times Q(k,y) \right)}}\\
    & = {} {\sum_{k=1}^{n}{\left(P(x,k) \times Q(k,y)\right)}} + {\sum_{k=1}^{n}{\left(S(x,k) \times Q(k,y)\right)}}\\
    & = {} P(x,y)\otimes Q(x,y)+S(x,y)\otimes Q(x,y). && 
\end{flalign*}
So the polynomials ${\left(P(x,y)+S(x,y) \right)}\otimes Q(x,y)$ and $P(x,y)\otimes Q(x,y)$+$S(x,y)$ $\otimes Q(x,y)$ are equal, entry-by-entry, and by the definition of polynomial equality, they are equal polynomials in $\mathcal{P}_{m}^{p}$. This completes the proof.

\begin{definition}[Identity under DP product]\label{definition:identity1}
Let $P \in \mathcal{P}_{m}^{n}$ be any polynomial. Then, the polynomial 
\begin{enumerate}[label=(\roman*)]
    \item $I_n \in \mathcal{P}_{n}^{n}$ is said to right identity, if $P(x,y)\otimes I_n(x,y)=P(x,y)$,
    \item $I_m \in \mathcal{P}_{m}^{m}$ is said to be left identity, if $I_m(x,y)\otimes P(x,y)=P(x,y)$.
\end{enumerate}
In general, the polynomial $I_n \in \mathcal{P}_{n}^{n}$ is said to be the identity, if 
\begin{equation}\label{definition:identity}
        P(x,y)\otimes I_n(x,y)=P(x,y)=I_n(x,y)\otimes P(x,y)
\end{equation}
for all $P \in \mathcal{P}_{n}^{n}$. 
\end{definition}

\begin{example}
Let $P\in \mathcal{P}_{3}^{2}$, $Q\in \mathcal{P}_{2}^{2}$, and $R \in \mathcal{P}_{3}^{3}$ and be three polynomials given by $P(x,y) = ax^2y+bxy+cx^2+dx+ey+f$, $Q(x,y)=2xy-3x-3y+5$, and $R(x,y) = \frac{1}{2}(3x^2y^2-12x^2y-12xy^2+10x^2+49xy+10y^2-42x-42y+38)$ respectively. Then, the definition \eqref{definition:dp product} implies that $P(x,y) \otimes Q(x,y)= ax^2y+bxy+cx^2+dx+ey+f$ and $R(x,y) \otimes P(x,y)= ax^2y+bxy+cx^2+dx+ey+f$.
%\begin{flalign*}
%    & P(x,y) \otimes Q(x,y)= ax^2y+bxy+cx^2+dx+ey+f, \\
%    \textrm{and } & R(x,y) \otimes P(x,y)= ax^2y+bxy+cx^2+dx+ey+f. &&
%\end{flalign*}
Therefore, $Q\in \mathcal{P}_{2}^{2}$ and $R \in \mathcal{P}_{3}^{3}$ are the right and left identity of $P\in \mathcal{P}_{3}^{2}$ respectively.
\end{example}

\begin{example}
Let $P(x,y) = axy+bx+cy+d$ and $Q(x,y)=2xy-3x-3y+5$ be two polynomials in $\mathcal{P}_{2}^{2}$. Then, the definition \eqref{definition:dp product} implies that
\begin{equation}
    P(x,y)\otimes Q(x,y)= axy+bx+cy+d = Q(x,y)\otimes P(x,y). \nonumber
\end{equation}
Thus, $P(x,y)\otimes Q(x,y)=P(x,y)=Q(x,y)\otimes P(x,y)$ holds. Hence, $Q(x,y)=2xy-3x-3y+5$ is the identity polynomial in $\mathcal{P}_{2}^{2}$.
\end{example}

\begin{remark}\label{identity polynomial remark}
The polynomials $I_1(x,y) = 1$, $I_2(x,y) = 2xy-3x-3y+5$, and $I_3(x,y) = \frac{1}{2}(3x^2y^2-12x^2y-12xy^2+10x^2+49xy+10y^2-42x-42y+38)$ are the identity polynomials in $\mathcal{P}_{1}^{1}$, $\mathcal{P}_{2}^{2}$ and $\mathcal{P}_{3}^{3}$ respectively.
\end{remark}

\begin{definition}[Inverse under DP product]\label{definition:inverse}
Let $P \in \mathcal{P}_{n}^{n}$ be any polynomial. Then, the polynomial $Q \in \mathcal{P}_{n}^{n}$ is said to be the inverse of $P(x,y)$, if
\begin{equation}
        P(x,y)\otimes Q(x,y)=I_n(x,y)=Q(x,y)\otimes P(x,y),
\end{equation}
where $I_n \in \mathcal{P}_{n}^{n}$ is the identity polynomial.
\end{definition}

\begin{remark}\label{definition:singular}
Let $\alpha_1,\alpha_2,\ldots,\alpha_n$ be $n$ scalars. The polynomial $P \in \mathcal{P}_{n}^{n}$ is invertible (or non-singular), if
\begin{equation}
        {\sum_{i=1}^{n}{\alpha_i.P(i,y)}}=0 \Rightarrow \alpha_i=0  \textrm{ for all } i=1,2,\ldots,n,
\end{equation}
\qquad \qquad \qquad \qquad \qquad \qquad \qquad or
\begin{equation}
        {\sum_{j=1}^{n}{\alpha_j.P(x,j)}}=0 \Rightarrow \alpha_j=0 \textrm{ for all } j=1,2,\ldots,m.
\end{equation}
\end{remark}

\begin{example}\label{ex:invertible}
Let $P(x,y)=-10xy+14x+13y-18$ and $Q(x,y)=15xy-21x-20y+28$ be two polynomials in $\mathcal{P}_{2}^{2}$. Suppose $\alpha_1,\alpha_2$ be two scalars, then ${\sum_{i=1}^{2}{\alpha_i.P(i,y)}}=0$ and ${\sum_{j=1}^{2}{\alpha_j.P(x,j)}}=0$ implies that $\alpha_1=\alpha_2=0$. Therefore, the polynomial $P \in \mathcal{P}_{2}^{2}$ is invertible. The, the polynomial $Q \in \mathcal{P}_{2}^{2}$ is not invertible, since ${\sum_{i=1}^{2}{\alpha_i.P(i,y)}}=0$ implies that $(\alpha_1,\alpha_2)=(2k,k)$, where $k$ is an arbitrary constant. Suppose $R \in \mathcal{P}_{2}^{2}$ be the inverse of $P \in \mathcal{P}_{2}^{2}$, given as
    \begin{equation}
     	R(x,y)=\beta_1 xy + \beta_2 x+\beta_3 y + \beta_4; \quad \beta_1,\beta_2,\beta_3,\beta_4 \in \mathbb{R}.
    \end{equation}
Since $I_2(x,y) = 2xy-3x-3y+5$ is the identity in $\mathcal{P}_{2}^{2}$, therefore the definition \eqref{definition:inverse} implies that $\beta_1 = 0$, $\beta_2 = -\frac{1}{2}$, $\beta_3 = -1$, and $\beta_4 = \frac{7}{2}$. Hence, $R(x,y)= -\frac{1}{2}(x+2y-7)$ is the inverse of the given polynomial $P \in \mathcal{P}_{2}^{2}$.
\end{example}

\begin{definition}[Eigenvalue of the polynomial]\label{definition:eigen value}
Let $P \in \mathcal{P}_{n}^{n}$ and $X \in \mathcal{P}_{n}^{1}$ be two polynomials. Suppose $\lambda$ is a scalar, if 
\begin{equation}
        P(x,y)\otimes X(x,1)= \lambda . X(x,1),
\end{equation}
then $\lambda$ is said to be ``eigenvalue" of the polynomial $P \in \mathcal{P}_{n}^{n}$ and $X \in \mathcal{P}_{n}^{1}$ is the ``eigen-polynomial" corresponding to the eigenvalue $\lambda$.
\end{definition}

\begin{example}
Let $P \in \mathcal{P}_{2}^{2}$ and $X_1,X_2 \in \mathcal{P}_{2}^{1}$ be three polynomials such that $P(x,y)=-xy+4x+3y-5$, $X_1(x,1)=5x-8$ and $X_2(x,1)=x$. Then, 
\begin{equation}
        P(x,y)\otimes X_1(x,1)= -X_1(x,1) \textrm{ and } P(x,y)\otimes X_2(x,1)= 7X_2(x,1). \nonumber
\end{equation}
Hence, $\lambda=-1$ and $\lambda=7$ are the eigenvalues of the polynomial $P(x,y)=-xy+4x+3y-5$, where $X_1(x,1)=5x-8$ and $X_2(x,1)=x$ are the eigen-polynomials corresponding to the eigenvalues $\lambda=-1$ and $\lambda=7$ respectively.
\end{example}

\begin{definition}[Power of the polynomial]\label{definition:power}
Let $P \in \mathcal{P}_{n}^{n}$ and $r \geq 2$, $r \in \mathbb{N}$, then
\begin{equation}
    {\left(P(x,y)\right)}^2= P(x,y)\otimes P(x,y) \in \mathcal{P}_{n}^{n}.
\end{equation}
\end{definition}

\begin{remark}
The polynomial $P \in \mathcal{P}_{n}^{n}$ is said to 
\begin{enumerate}
    \item involuntary, if ${\left(P(x,y)\right)}^2= I_n(x,y)$, where $I_n(x,y)$ is the identity polynomial in $\mathcal{P}_{n}^{n}$.
    \item idempotent, if ${\left(P(x,y)\right)}^2= P(x,y)$.
    \item nilpotent with index $r$, if ${\left(P(x,y)\right)}^r = 0$, where $r$ is the least positive integer.
    \item periodic with index $r$, if ${\left(P(x,y)\right)}^{r+1} = P(x,y)$, where $r$ is the least positive integer.
\end{enumerate}
\end{remark}

\subsection{Transpose, symmetry, skew-symmetry, and orthogonality of the polynomials in $\mathcal{P}_{m}^{n}$}
\begin{definition}[Transpose of the polynomial]\label{definition:transpose}
Let $P(x,y)$ be a polynomial in $\mathcal{P}_{m}^{n}$, then its transpose is defined as 
\begin{equation}
    P^T(x,y)= P(y,x), 
\end{equation}
such that $P^T \in \mathcal{P}_{n}^{m}$.
\end{definition}

\begin{remark}
The polynomial $P \in \mathcal{P}_{n}^{n}$ is said to 
\begin{enumerate}
    \item symmetric, if $P(x,y)=P(y,x)$.
    \item skew-symmetric, if $P(x,y) = - P(y,x)$.
    \item orthogonal, if $P(x,y)\otimes P^{T}(x,y)= I_n(x,y) = P^{T}(x,y)\otimes P(x,y)$, where $I_n(x,y)$ is the identity polynomial in $\mathcal{P}_{n}^{n}$.
\end{enumerate}
\end{remark}

%Some particular cases can be listed as:
%\begin{itemize}
%\item 
%$\P_{1}^{1}:= \left \{
%{\lambda _{00}} \mid {\lambda _{00}}\in \mathbb{R}\right \}$
%
%\item 
%$\P_{1}^{2}:= \left \{{\lambda _{00}}+{\lambda _{01}y} \left|  %\lambda _{00}, \lambda _{01}\in \mathbb{R}
%\right. \right \}$
%
%\item 
%$\P_{2}^{1}:= \left \{{\lambda _{00}}+{\lambda _{10}x} \left| \lambda _{00}, \lambda _{10} \in \mathbb{R} 
%\right. \right \}$
%
%\item 
%$\P_{2}^{2}:= \left \{
%\begin{array}{ll}
%    {\lambda _{00}+\lambda _{01}y + \lambda _{10}x + \lambda _{11}xy}  
%\end{array}
%\left| 
%\begin{array}{ll}
%    {\lambda _{{\alpha}{\beta}}}\in \mathbb{R} \textrm{ for all } \alpha,\beta \in \left \{0,1 \right \}
%\end{array}
%\right. \right \}$
%
%\item $\P_{2}^{3}:= \left \{\lambda _{00}+\lambda _{01}y+\lambda _{02}y^2+\lambda _{10}x+\lambda _{11}xy+\lambda _{12}{x}{y^2}
%\left| \begin{array}{ll} {\lambda _{{\alpha}{\beta}}}\in \mathbb{R} \textrm{ for all}\\
%\textrm{ } \alpha \in \left \{0,1 \right \},\\
%   \beta \in \left \{0,1,2\right \}
%\end{array}
%\right.
%\right \}$
%\end{itemize}
%

\section{Existence of the subspace $\mathcal{P}_{m}^{n}$, poisedness of the MIP, and some polynomial construction formulae}
In this section, the existence of the subspace $\mathcal{P}_{m}^{n}$ against the solution of the MIP \eqref{mip} and the poisedness of the MIP \eqref{mip} in it will be proved. Three polynomial construction formulae will be derived using univariate Lagrange, Newton-forward difference, and Newton-backward difference interpolation formulae \cite{Atkinson:1,Jain:1} with respect to the MIP\eqref{mip} for a given matrix $(a_{ij}) \in \mathbb{R}^{m \times n}$ in the subspace $\mathcal{P}_{m}^{n}$. Some arguments and other construction methods will be remarked. 

\begin{theorem} \label{thm:dp} 
Let $a_{ij} \in \mathbb{R}^{m \times n}$ be any given matrix. Then, there exists a unique polynomial $P \in \mathcal{P}_{m}^{n}$ with respect to the set $\mathcal{D}$ which satisfy the MIP \eqref{mip} .
\end{theorem}

\noindent \textbf{Proof:} The proof consist of two parts, existence and uniqueness. Both will be proved one by one as follows:

\noindent \textit{Existence: } Let $A=(a_{ij}) \in \mathbb{R}^{m \times n}$ be a given matrix. Then, for the given set of nodes 
	\begin{equation}\label{rthrow:nodes}
            \Omega_m = \left\{i:i=1,2,\ldots,n\right\},
	\end{equation}
associated with the $r$th column, there exists a unique polynomial $p_r(x)$ in one variable $x$, of degree at most $m-1$, for all $r=1,2,\ldots,m$, satisfying the problem
	\begin{equation}\label{intp:1}
            p_r(i)=a_{ir} \textrm{ for all } i=1,2,\ldots,m.
	\end{equation}
Therefore, the matrix $A=(a_{ij})_{m\times n}$ can be represented  as $\begin{pmatrix}
        p_1(x) p_2(x) \ldots p_m(x)
        \end{pmatrix} _{1 \times n}$. Again, for the given set of nodes 
	\begin{equation}\label{gen:nodes}
            \Omega_n = \left\{j:j=1,2,\ldots,n\right\} ,  
	\end{equation}
there exists a unique polynomial $P(x,y)$ in the variable $y$ (for fixed $x$), of degree at most $n-1$, satisfying the problem
	\begin{equation}\label{intp:2}
            P(x,j)=p_j(x) \textrm{ for all } j=1,2,\ldots,n.
	\end{equation}
On combing, $P(x,y)$ is a polynomial in $\mathcal{P}_{m}^{n}$, satisfying the MIP (\ref{mip}) with respect to the set $\mathcal{D} =\Omega_m \times \Omega_n$. 
%It should be noted that the degree of the variables $x$ and $y$ is at most $m-1$ and $n-1$ respectively. 
This completes the proof of existence part.

\noindent \textit{Uniqueness:} The polynomial expressions with respect to the set of nodes \eqref{rthrow:nodes} can be written as 
\begin{equation}
        a_{ir}=p_r(i)+R_{m,r}(i), \;\;i=1,2,\ldots,m,
\end{equation}
where $R_{m,r}(i)$ are the remainders for all $r=1,2,\ldots,n$. The polynomials, satisfy the conditions $a_{kr}=p_r(k)$ for all $k=1, 2, \ldots, m$. Therefore, the remainder term $R_{m,r}$; $r=1,2,\ldots,n$ vanishes at $i=k$ for all $k=1, 2, \ldots, m$. Thus, the polynomials satisfying the conditions (\ref{intp:1}) are unique for all $r=1,2,\ldots,n$. Again, the polynomial expression with respect to the set of nodes (\ref{gen:nodes}) can be written as 
\begin{equation}
        p_j(x)=P(x,j)+R_n(j),\;\;j=1,2,\ldots,n,
\end{equation}
where $R_n(j)$ is the remainder. The polynomial $P(x,j)$, satisfies the conditions $p_r(k)=P(x,k)$ for all $k=1, 2, \ldots, n$. Therefore, the remainder term $R_n$ vanishes at $j=k$ for all $k=1, 2, \ldots, n$. Hence, polynomial satisfying the conditions (\ref{intp:2}) is unique. This completes the theorem.

\begin{remark}
We have used an approach of uni-variate polynomial interpolation for the considered MIP\eqref{mip} and the interpolation space coincide with the space suggested by the author Narumi in the year 1920 in \cite{narumi1920some}. Narumi had not mentioned the used approach of interpolation for the existence of the space, however a method of bivariate divided difference is used to evaluate the coefficients of the given polynomial form in the suggested space.
\end{remark}

\begin{remark}
In general, the interpolation methods assume some certain properties on the sought structure such as smoothness, regularity (nonexistence of some singularities) in a sufficiently large vicinity and so on. The ignorance of these types of assumptions would be dangerous from the error point of view, \cite{mft1,mft2,mft3}. We are assuming that the original function is smooth and regular or `well-behaved' in the given domain. The present work is carried out in the algebraic direction and approximation or error estimation is not the part of the article.
\end{remark}

%\begin{corollary}
%For all $A=(a_{ij}) \in \mathbb{R}^{m \times n}$, there exists a unique $p_r \in \Pi_{n-1}^1$, associated with the $r$th row, $r=1,2,\ldots,m$ of the matrix, satisfying $a_{rj}=p_r(j)$ for all $j=1,2,\ldots,n$.
%\end{corollary}

%\begin{corollary}
%For all $A=(a_{ij}) \in \mathbb{R}^{m \times n}$, there exists a unique $p_s \in \Pi_{m-1}^1$, associated the $s$th column, $s=1,2,\ldots,n$ of the matrix $A$, satisfying $a_{is}=p_s(i)$ for all $i=1,2,\ldots,m$.
%\end{corollary}

In the next theorems, some formulae to construct the unique polynomial $P \in \mathcal{P}_{m}^{n}$ satisfying the MIP \eqref{mip} with respect to the set $\mathcal{D}$ for a given matrix in $\mathbb{R}^{m \times n}$ are formulated using well-known univariate interpolation formulae.

\begin{theorem} \label{thm:construction by lagrange} 
Let $A=(a_{ij}) \in \mathbb{R}^{m \times n}$ be a given matrix, then there exist a unique polynomial $P \in \mathcal{P}_{m}^{n}$, satisfying $a_{ij}=P(i,j)$ for all $i=1,2,\ldots,m$, $j=1,2,\ldots,n$, given by
	\begin{equation} \label{lag:mat}
		P(x,y)=\sum_{r=1}^{n} \left \{ \prod_{\alpha=1, \alpha \neq r}^{n} {\frac {(y-\alpha)}{(r-\alpha)}} \right \}{p_r(x)},
	\end{equation}
where
	\begin{equation} \label{lag:row}
        \quad p_r(x)=\sum_{k=1}^{m} \left \{ \prod_{\alpha=1, \alpha \neq k}^{m} {\frac {(x-\alpha)}{(k-\alpha)}} \right \}{a_{kr}} \textrm{ for all } r=1,2,\ldots,n.
	\end{equation}
\end{theorem}

\noindent \textbf{Proof:} Let $A=(a_{ij}) \in \mathbb{R}^{m \times n}$ be a given matrix, then the poisedness of the MIP \eqref{mip} with respect to the set $\mathcal{D}$ of nodes insures that, there exist unique polynomial $P \in \mathcal{P}_{m}^{n}$ satisfying $a_{ij}=P(i,j)$ for all $i=1,2,\ldots,m$,\; $j=1,2,\ldots,n$. Using univariate Lagrange interpolation formula, the polynomials with respect to the set of nodes (\ref{rthrow:nodes}), satisfying the problem (\ref{intp:1}), can be written in the form 
\begin{equation} \label{lag:1}
        p_r(x)=l_{1,r}(x)a_{1r}+l_{2,r}(x)a_{2r}+\ldots+l_{n,r}(x)a_{nr},
\end{equation}
where $l_{k,r}(x)$, $k=1,2,\ldots,m$ are the polynomials of one variable $x$, of degree at most $m-1$ for all $r=1,2,\ldots,n$. The polynomials (\ref{lag:1}) will satisfy the interpolating conditions (\ref{intp:1}), if and only if 
\begin{equation} \label{lag:2}
        l_{k,r}(x_i)=\begin{cases}
                    0,  & \text{if } k \neq i\\
                    1,  & \text{if } k =i 
                \end{cases} \textrm{ for all } r=1,2,\ldots,n,
\end{equation}
and the polynomials $l_{k,r}(x)$ for all $r=1,2,\ldots,n$, satisfying the conditions (\ref{lag:2}), can be given as
\begin{equation} \label{lag:3}
        l_{k,r}(x)= \prod_{\alpha=1, \alpha \neq k}^{m} {\frac {(x-\alpha)}{(k-\alpha)}} \textrm{ for all } k=1,2,\ldots,m.
\end{equation}
The equations (\ref{lag:1}) and (\ref{lag:3}), insured that
	\begin{equation} \label{lag:row2}
        \quad p_r(x)=\sum_{k=1}^{m} \left \{ \prod_{\alpha=1, \alpha \neq k}^{m} {\frac {(x-\alpha)}{(k-\alpha)}} \right \}{a_{kr}} \textrm{ for all } r=1,2,\ldots,n. \nonumber
	\end{equation}
Again, using univariate Lagrange interpolation formula, the polynomial with respect to the set nodes (\ref{gen:nodes}), satisfying the problem (\ref{intp:2}), can be written in the form 
\begin{equation} \label{lag:4}
        P(x,y)=L_1(y)p_1(x)+L_2(y)p_2(x)+\ldots+L_n(y)p_n(x),
\end{equation}
where $L_r(x)$, $r=1,2,\ldots,n$ are the polynomials of one variables $y$, of degree at most $n-1$, i.e., $L_r \in \Pi_{n-1}^1$. The polynomial (\ref{lag:4}) will satisfy the interpolating conditions (\ref{intp:2}), if and only if 
\begin{equation} \label{lag:5}
        L_r(y_j)=\begin{cases}
                    0,  & \text{if } r \neq j\\
                    1,  & \text{if } r = j 
                \end{cases},
\end{equation}
and the polynomials $L_r(y)$, satisfying the conditions (\ref{lag:5}), can be given as
\begin{equation} \label{lag:6}
        L_r(y)= \prod_{\alpha=1, \alpha \neq r}^{n} {\frac {(y-\alpha)}{(r-\alpha)}} \textrm{ for all } r=1,2,\ldots,n.
\end{equation}
The combination of the equations \eqref{lag:4} and \eqref{lag:6}, completes the result.
%The equations (\ref{lag:4}) and (\ref{lag:6}), insured that 
%	\begin{equation} \label{lag:mat2}
%		P(x,y)=\sum_{r=1}^{m} \left \{ \prod_{\alpha=1, \alpha \neq r}^{m} {\frac {(x-\alpha)}{(r-\alpha)}} \right \}{p_r(y)}.
%	\end{equation}

In the similar manner, two more formulae using univariate Newton forward and backward interpolation formulae, for the given set of nodes \eqref{rthrow:nodes} and \eqref{gen:nodes} satisfying the interpolation problems \eqref{intp:1} and \eqref{intp:2} respectively, can be given as follows: 

\begin{theorem} \label{thm:construction by newton forwatd} 
Let $A=(a_{ij}) \in \mathbb{R}^{m \times n}$ be a given matrix, then there exist a unique polynomial $P \in \mathcal{P}_{m}^{n}$, satisfying $a_{ij}=P(i,j)$ for all $i=1,2,\ldots,m$, $j=1,2,\ldots,n$, given by
\begin{equation} \label{gnf:3}
\begin{aligned}
        P(x,y) = {} & p_1(x) + \frac{(y-1)}{1!}{\Delta p_1(x)} + \frac{(y-1)(y-2)}{2!}{\Delta^2 p_1(x)}\\
        & \quad \quad \quad \quad +\ldots+\frac{(y-1)(y-2)\ldots(y-(n-1))}{(n-1)!}{\Delta^{n-1} p_1(x)},
\end{aligned}
\end{equation}
where
\begin{equation} \label{gnf:4}
    \Delta^{r} p_1(x)=\sum_{\beta=0}^{r}{(-1)^{\beta}}{r \choose \beta}{p_{r+1-\beta}(x)} \textrm{ for all } r=1,2,\ldots,n-1,
\end{equation}
\newline
and
\begin{equation} \label{gnf:1}
\begin{aligned}
        p_r(x) = {} & a_{1r} + \frac{(x-1)}{1!}{\Delta a_{1r}} + \frac{(x-1)(x-2)}{2!}{\Delta^2 a_{1r}}\\
        & \quad \quad \quad \quad \quad \quad +\ldots+\frac{(x-1)(x-2)\ldots(x-(m-1))}{(m-1)!}{\Delta^{m-1} a_{1r}},
\end{aligned}
\end{equation}
for all $r=1,2,\ldots,n$, where
\begin{equation} \label{gnf:2}
    \Delta^{k} a_{1r}=\sum_{\alpha=0}^{k}{(-1)^{\alpha}}{k \choose \alpha}{a_{(k+1-\alpha)r}} \textrm{ for all } k=1,2,\ldots,m-1.
\end{equation}
\end{theorem}

\begin{theorem} \label{thm:construction by newton backward} 
Let $A=(a_{ij}) \in \mathbb{R}^{m \times n}$ be a given matrix, then there exist a unique polynomial $P \in \mathcal{P}_{m}^{n}$, satisfying $a_{ij}=P(i,j)$ for all $i=1,2,\ldots,m$, $j=1,2,\ldots,n$, given by
\begin{equation} \label{gnb:3}
\begin{aligned}
        P(x,y) = {} & p_n(x) + \frac{(y-n)}{1!}{\nabla p_n(x)} + \frac{(y-n)(y-(n-1))}{2!}{\nabla^2 p_n(x)}\\
        & \quad \quad \quad \quad +\ldots+\frac{(y-n)(y-(n-1))\ldots(y-2)}{(n-1)!}{\nabla^{n-1} p_n(x)},
\end{aligned}
\end{equation}
where
\begin{equation} \label{gnb:4}
    \nabla^{r} p_n(x)=\sum_{\beta=0}^{r}{(-1)^{\beta}}{r \choose \beta}{p_{n-\beta}(x)} \textrm{ for all } r=1,2,\ldots,n-1.
\end{equation}
and
\begin{equation} \label{gnb:1}
\begin{aligned}
        p_r(x) = {} & a_{mr} + \frac{(x-m)}{1!}{\nabla a_{mr}} + \frac{(x-m)(x-(m-1))}{2!}{\nabla^2 a_{nr}}\\
        & \quad \quad \quad \quad \quad \quad +\ldots+\frac{(x-m)(x-(m-1))\ldots(x-2)}{(m-1)!}{\nabla^{m-1} a_{mr}},
\end{aligned}
\end{equation}
for all $r=1,2,\ldots,n$, where
\begin{equation} \label{gnb:2}
    \nabla^{k} a_{mr}=\sum_{\alpha=0}^{k}{(-1)^{\alpha}}{k \choose \alpha}{a_{(m-\alpha)r}} \textrm{ for all } k=1,2,\ldots,m-1.
\end{equation}
\end{theorem}

\begin{remark}
Suppose $P \in \mathcal{P}_{m}^{n}$ is the unique polynomial of the form \eqref{dpsosk} which satisfy the MIP \eqref{mip} for the given matrix $(a_{ij}) \in \mathbb{R}^{m \times n}$. Then, the $mn$ independent conditions $a_{ij}=P(i,j)$ for all $i=1,2,\ldots,m$,\; $j=1,2,\ldots,n$ constitute a non-homogeneous linear system of $mn$ equations of the form ${\Lambda}X=\mu$, where $X=(\lambda_{00},\lambda_{01},\ldots,\lambda_{(m-1)(n-1)})^T$, $\mu = (a_{11},a_{12},\ldots,a_{mn})^T$, and $\Lambda$ is a coefficient matrix. The uniqueness of the solution implies that $det(\Lambda) \neq 0$. Using some suitable method to solve linear system of non-homogeneous equations, the coefficients $\lambda_{k_1k_2}$, where $k_1=0,1,\ldots,m-1$ and $k_2=0,1,\ldots,n-1$ can be determined uniquely. 
\end{remark}

\section{Isomorphism and isomorphic structures between the spaces $\mathbb{R}^{m \times n}$ and $\mathcal{P}_{m}^{n}$}
Suppose the notation $P_A(x,y)$ represent the unique polynomial $P_A \in \mathcal{P}_{m}^{n}$ which satisfy the MIP \eqref{mip} for the given matrix $A \in \mathbb{R}^{m \times n}$.
%constructed by any one of the formula given in the theorems \ref{thm::construction by lagrange}, \ref{thm::construction by newton forwatd}, and \ref{thm::construction by newton backward} or any other method. 
This notation will be used frequently in the reaming part of the paper. The spaces $\mathbb{R}^{m \times n}$ and $\mathcal{P}_{m}^{n}$ are isomorphic with respect to the vector spaces structure, therefore there must exist an isomorphism (invertible linear map) between them \cite{axler15}. 

Using the result of the theorem \ref{thm:dp}, the polynomial map $D_p:\mathbb{R}^{m \times n} \rightarrow \mathcal{P}_{m}^{n}$ can be defined as 
\begin{equation}\label{dp:map}
D_p(A) = P_A(x,y) \textrm{ for all } A \in \mathbb{R}^{m \times n}. 
\end{equation}
Again, let the map $f:\mathbb{R}^{m \times n} \times \mathbb{R}^{n \times q} \rightarrow \mathbb{R}^{m \times q}$ is defined as $f(A,B)=A.B$ for all $A \in \mathbb{R}^{m \times n}$, $B \in \mathbb{R}^{n \times q}$, then the composition map $D_p o f:\mathbb{R}^{m \times n} \times \mathbb{R}^{n \times q} \rightarrow \mathbb{R}^{m \times q}$
is given as $D_p o f(A,B)=P_{A.B}(x,y)$ for all $A \in \mathbb{R}^{m \times n}$, $B \in \mathbb{R}^{n \times q}$. Also, let the maps $g:\mathbb{R}^{m \times n} \times \mathbb{R}^{n \times q} \rightarrow \mathcal{P}_{m}^{n} \times \mathcal{P}_{n}^{q}$ and $h:\mathcal{P}_{m}^{n} \times \mathcal{P}_{n}^{q} \rightarrow \mathcal{P}_{m}^{q}$ are defined as $g(A,B)=(P_A(x,y),P_B(x,y))$ for all $A \in \mathbb{R}^{m \times n}$, $B \in \mathbb{R}^{n \times q}$ and  $h(P,Q) = P(x,y) \otimes Q(x,y)$ for all $P \in \mathcal{P}_{m}^{n}$, $Q \in \mathcal{P}_{n}^{q}$ respectively. Then, the composition map $h o g:\mathbb{R}^{m \times n} \times \mathbb{R}^{n \times q} \rightarrow \mathcal{P}_{m}^{q}$ is given as $h o g(A,B)=P_A(x,y) \otimes P_B(x,y)$ for all $A \in \mathbb{R}^{m \times n}$, $B \in \mathbb{R}^{n \times q}$. 

In this section, by proving the invertibility and linearity, it will be proved that the polynomial map \eqref{dp:map} is an isomorphism. Further, the property $P_{A.B}(x,y)= P_{A}(x,y)\otimes P_{B}(x,y)$ for all $A \in \mathbb{R}^{m \times n}$, $B \in \mathbb{R}^{n \times q}$ will be proved to establish the discussed product structure. At last, the existence of the identity polynomial under the operation `$\otimes$' will be ensured and it will be proved that the algebraic structures $\langle \mathbb{R}^{n \times n}, .\rangle$ and $\langle \mathcal{P}_{n}^{n}, \otimes  \rangle$ are isomorphic.

\begin{theorem}\label{thm:inversedp}
       For all $P \in \mathcal{P}_{m}^{n}$, there exists a unique matrix $A =(a_{ij}) \in \mathbb{R}^{m \times n}$ such that $a_{ij}=P(i,j)$ for all $i=1,2,\ldots,m$, $j=1,2,\ldots,n$.
\end{theorem}

\begin{theorem} \label{thm:linearity}
Let $A,B \in \mathbb{R}^{m \times n}$ and $\alpha$ be a real scalar, then the following properties holds:
\begin{enumerate}[label=(\roman*)]
    \item $P_{A+B}(x,y)= P_{A}(x,y)+P_B(x,y)$ in $\mathcal{P}_{m}^{n}$.
    \item $P_{\alpha A}(x,y)= \alpha P_{A}(x,y)$ in $\mathcal{P}_{m}^{n}$.
\end{enumerate}
\end{theorem}
\noindent \textbf{Proof:} 
Let $A,B \in \mathbb{R}^{m \times n}$ be two matrices given as $A=(a_{ij})_{m \times n}$ and $B =(b_{ij})_{m \times n}$. Then, $A+B = (a_{ij}+b_{ij})_{m \times n}$ and $\alpha A = (\alpha a_{ij})_{m \times n}$, for some scalar $\alpha$. 
Therefore, there exists unique $P_A \in \mathcal{P}_{m}^{n}$, $P_B \in \mathcal{P}_{m}^{n}$, $P_{A+B} \in \mathcal{P}_{m}^{n}$, and $P_{\alpha A} \in \mathcal{P}_{m}^{n}$, which satisfy the MIP's
\begin{equation}
P_A(i,j)=a_{ij} \textrm{ for all } i=1,2,\ldots,m,\; j=1,2,\ldots,n,
\end{equation}
\begin{equation}
    P_B(i,j)=b_{ij} \textrm{ for all } i=1,2,\ldots,m,\; j=1,2,\ldots,n,
\end{equation}
\begin{equation}\label{mip:addition}
P_{A+B}(i,j)= a_{ij} + b_{ij} \textrm{ for all } i=1,2,\ldots,m,\; j=1,2,\ldots,n,
\end{equation}
\begin{equation}\label{mip:scalar multiplication}
\textrm{and } P_{\alpha A}(i,j)= \alpha a_{ij} \textrm{ for all } i=1,2,\ldots,m,\; j=1,2,\ldots,n,
\end{equation}
respectively. 
%$\{(i,j):i=1,2,\ldots,m,$ $j=1,2,\ldots,n\}$.
Again, the definition of usual addition %\eqref{definition:addition} 
and scalar multiplication 
%\eqref{definition:scalar multiplication} 
of the polynomials in $\mathcal{P}_{m}^{n}$ implies that 
\begin{flalign*}
(i) \;\; P_A(i,j)+P_B(i,j) & = a_{ij}+b_{ij} \textrm{ for all } i=1,2,\ldots,m,\; j=1,2,\ldots,n \\
&= P_{a+B}(i,j) \textrm{ for all } i=1,2,\ldots,m,\; j=1,2,\ldots,n, &&
\end{flalign*}
\begin{flalign*}
and \;\; (ii) \;\; \alpha P_{A}(i,j) & = \alpha a_{ij} \textrm{ for all } i=1,2,\ldots,m,\; j=1,2,\ldots,n\\
& = \alpha P_{\alpha A}(i,j) \textrm{ for all } i=1,2,\ldots,m,\; j=1,2,\ldots,n.&&
\end{flalign*}
The MIP \eqref{mip:addition} and \eqref{mip:scalar multiplication} have unique solution in the space $\mathcal{P}_{m}^{n}$ with respect to the set $\mathcal{D}$ by the theorem \ref{thm:dp}. This completes the proof.

\begin{remark}
If $A_1,A_2,...,A_k \in M_{n,n}$, and $\alpha_1,\alpha_2,...,\alpha_k$ be some real scalars, then 
\begin{equation}
    P_{\alpha_1{A_1}+\alpha_2{A_2}+...+\alpha_k{A_k}}(x,y)= {\sum_{i=1}^{k}}{\alpha_i {P_{A_i}(x,y)}}. \nonumber
\end{equation}
\end{remark}

On combining the theorems \ref{thm:dp}, \ref{thm:inversedp}, and \ref{thm:linearity}, the properties of the polynomial map \eqref{dp:map} can be combine as follows.

\begin{theorem}\label{thm:isomorphic map}
The polynomial map $D_p:\mathbb{R}^{m \times n} \rightarrow \mathcal{P}_{m}^{n}$ defined by 
\begin{equation}\label{linearity}
	   D_p(A) = P_A(x,y) \textrm{ for all } A \in \mathbb{R}^{m \times n}.
\end{equation}
is an ispmprphism.
\end{theorem}

\begin{remark}
The inverse linear map ${D_p}^{-1}:\mathcal{P}_{m}^{n} \rightarrow \mathbb{R}^{m \times n}$ is given by
\begin{equation}
	   {D_p}^{-1}(P(x,y)) = {\left[P(i,j)\right]}_{m \times n} \textrm{ for all } P \in \mathcal{P}_{m}^{n}.
\end{equation}
\end{remark}

\begin{theorem}\label{thm:product structure}
Let $A \in \mathbb{R}^{m \times n}$ and $B \in \mathbb{R}^{n \times q}$, then
$P_{AB}(x,y)= P_{A}(x,y)\otimes P_{B}(x,y) \textrm{ in } \mathcal{P}_{m}^{q}$.
\end{theorem}
\noindent \textbf{Proof:} Let $A \in \mathbb{R}^{m \times n}$ and $B \in \mathbb{R}^{n \times q}$ be two matrices given as $A=(a_{ij})_{m \times n}$ and $B =(b_{ij})_{n \times p}$. Then, 
\begin{equation}
        AB = \left({\sum_{k=1}^{n}}{\left(a_{ik} \times b_{kj} \right)} \right)_{m \times q} \nonumber
\end{equation}
Therefore, there exists unique $P_A \in \mathcal{P}_{m}^{n}$, $P_B \in \mathcal{P}_{m}^{n}$ and $P_{AB} \in \mathcal{P}_{m}^{p}$ satisfying the MIP's
\begin{equation}
P_A(i,j)=a_{ij} \textrm{ for all } i=1,2,\ldots,m,\; j=1,2,\ldots,n,
\end{equation}
\begin{equation}
    P_B(i,j)=b_{ij} \textrm{ for all } i=1,2,\ldots,n,\; j=1,2,\ldots,q,
\end{equation}
\begin{equation}\label{mip:multiplication}
P_{AB}(i,j)= {\sum_{k=1}^{n}}{\left(a_{ik} \times b_{kj} \right)} \textrm{ for all } i=1,2,\ldots,m,\; j=1,2,\ldots,q,
\end{equation}
with respect to the sets $\{(i,j):i=1,2,\ldots,m$,\; $j=1,2,\ldots,n\}$, $\{(i,j):i=1,2,\ldots,n$,\; $j=1,2,\ldots,q\}$, and $\{(i,j):i=1,2,\ldots,m$,\; $j=1,2,\ldots,q\}$ respectively.
Again, the definition of DP product \eqref{definition:dp product} implies that 
\begin{flalign*}
P_A(i,j)\otimes P_B(i,j) & = {\sum_{k=1}^{n}}{\left(P_A(i,k) \times P_B(k,j) \right)} \textrm{ for all } i=1,2,\ldots,m,\; j=1,2,\ldots,p \\
&={\sum_{k=1}^{n}}{\left(a_{ik} \times b_{kj} \right)} \textrm{ for all } i=1,2,\ldots,m,\; j=1,2,\ldots,p\\
& = P_{AB}(i,j) \textrm{ for all } i=1,2,\ldots,m,\; j=1,2,\ldots,q.&&
\end{flalign*}
Using Theorem \ref{thm:dp}, the MIP \eqref{mip:multiplication} with respect to the set $\{(i,j):i=1,2,\ldots,m$,\; $j=1,2,\ldots,q\}$ will have unique solution in $\mathcal{P}_{m}^{q}$. This completes the proof.

\begin{remark}
If $A_1,A_2,...,A_k \in \mathbb{R}^{m \times n}$, $B_1,B_2,...,B_k \in M_{n,p}$, and $\alpha_1,\alpha_2,...,\alpha_k$ be some real scalars, then 
\begin{equation}
    P_{\alpha_1{A_1}{B_1}+\alpha_2{A_2}{B_2}+...+\alpha_k{A_k}{B_k}}(x,y)= {\sum_{i=1}^{k}}{\alpha_i {\left(P_{A_i}(x,y)\otimes P_{B_i}(x,y)\right)}}\in \mathcal{P}_{m}^{p}. \nonumber
\end{equation}
\end{remark}

\begin{corollary}[Existence of identity polynomial under DP product]\label{coro:identity polynomial}
        Suppose $P_A \in \mathcal{P}_{m}^{n}$, then
\begin{enumerate}[label=(\roman*)]
    \item there exist a $P_{I_n} \in \mathcal{P}_{n}^{n}$, such that $P_{A}(x,y)\otimes P_{I_n}(x,y)=P_{A}(x,y)$,
    \item there exist a $P_{I_m} \in \mathcal{P}_{m}^{m}$, such that $P_{I_m}(x,y)\otimes P_{A}(x,y)=P_{A}(x,y)$.
\end{enumerate}
\end{corollary}

The combination of the theorems \ref{thm:isomorphic map} and \ref{thm:product structure} establish a very important result, that can be summarize as follows.

\begin{theorem}\label{thm:isomorphic algebra structure}
The binary structures $\langle \mathbb{R}^{n \times n}, .\rangle$ and $\langle \mathcal{P}_{n}^{n}, \otimes  \rangle$ are isomorphic. 
\end{theorem}

For $m=n$, the combination of the theorems  \ref{thm:vector space}, \ref{thm:multiplicative closure}, \ref{thm:multiplicative associativity}, \ref{thm:distributivity across addition} and the corollary \ref{coro:identity polynomial} can be resulted as follows.
\begin{theorem}
The binary structure $\langle \mathcal{P}_{n}^{n},+,\otimes  \rangle$ form a ring with unity.
\end{theorem}

\section{Numerical verification}\label{section 2.1}
In this section, five examples will be included to illustrate and verify the results.

\begin{example}\label{ex:1}
Let $\delta \in \mathbb{R}^{1 \times 3}$, $\eta \in \mathbb{R}^{3 \times 1}$, $\vartheta, \psi \in \mathbb{R}^{2 \times 2}$, $\zeta \in \mathbb{R}^{2 \times 3}$ and $\theta \in \mathbb{R}^{3 \times 2}$ be six matrices 			
$\delta = \begin{pmatrix}			% THIS MATRIX GIVES NO SOLUTION IN LAGRANGE METHOD
		1 & -1 & -2
		\end{pmatrix}$,   
		$\eta = \begin{pmatrix}				% THIS MATRIX GIVES INFINITE SOLUTION IN LAGRANGE METHOD
		-1\\
		1\\
		3
		\end{pmatrix}$, 
		$\vartheta = \begin{pmatrix}     % LAGRANGE METHOD IS NOT APPLICABLE
		-15 & 36\\
		-1 & 96
		\end{pmatrix}$,
		$\psi = \begin{pmatrix}			% LAGRANGE METHOD IS NOT APPLICABLE
		-216 & 2\\
		540 & -5
		    \end{pmatrix}$,
		$\zeta = \begin{pmatrix}			% THIS MATRIX GIVES NO SOLUTION IN LAGRANGE METHOD
		1 & 0 & 2\\
		-1 & 2 & -3
		    \end{pmatrix}$ and 
		$\theta = \begin{pmatrix}		% THIS MATRIX GIVES INFINITELY MANY SOLUTIONS IN LAGRANGE METHOD
		1 & 0\\
		0 & -1\\
		1 & 0
		    \end{pmatrix}$ respectively. 
Therefore, there exists unique $P_{\delta}(1,y) \in \mathcal{P}_2^1$, $P_{\eta}(x,1) \in \mathcal{P}_2^1$, $P_{\vartheta}(x,y) \in \mathcal{P}_2^2$, $P_{\psi}(x,y) \in \mathcal{P}_2^2$, $P_{\zeta}(x,y) \in \mathcal{P}_2^3$ and $P_{\theta}(x,y) \in \mathcal{P}_3^2$ given by $P_{\delta}(1,y) = \frac{1}{2}y^2-\frac{7}{2}y+4$, $P_{\eta}(x,1) = 2x-3$, $P_{\vartheta}(x,y) = 46xy-32x+5y-34$, $P_{\psi}(x,y) = -763xy+1519x+981y-1953$, $P_{\zeta}(x,y) = -\frac{11}{2}xy^2+7y^2+\frac{41}{2}xy-17x-26y+22$, and $P_{\theta}(x,y) = x^2-4x-y+5$ respectively.

%Surface diagrams of the polynomials indicating the data points for the MIP's with respect to the matrices $\delta$, $\eta$, $\vartheta$, $\psi$, $\zeta$ and $\theta$ are represented in the figure \ref{matrix:6}.

\begin{figure}[H]
    \centering
    \begin{subfigure}[b]{0.3\textwidth}
        \includegraphics[width=\textwidth]{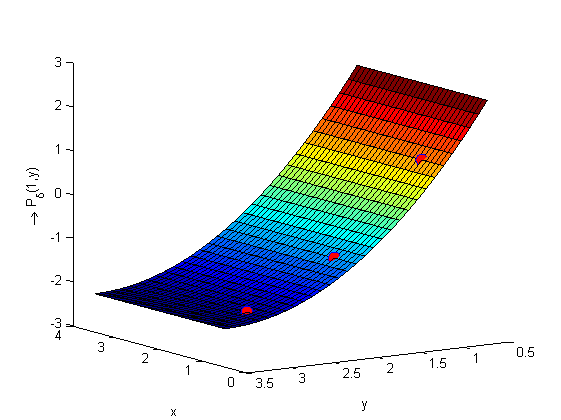}
        \caption{$P_{\delta} \in {\mathcal{P}_{2}^{1}}$}
        \label{fig:delta}
    \end{subfigure}
    ~ %add desired spacing between images, e. g. ~, \quad, \qquad, \hfill etc. 
    \begin{subfigure}[b]{0.3\textwidth}
        \includegraphics[width=\textwidth]{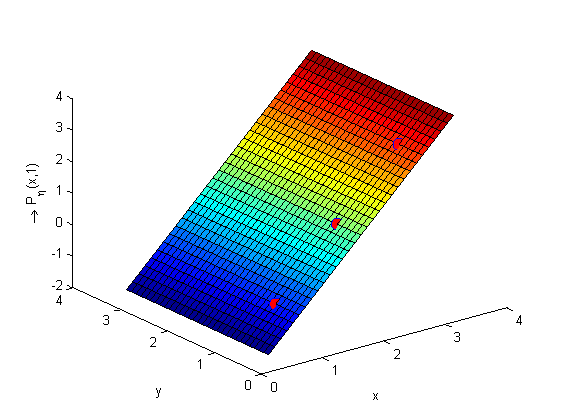}
        \caption{$P_{\eta} \in {\mathcal{P}_{2}^{1}}$}
        \label{fig:eta}
    \end{subfigure}
    %\caption{Graphs of the polynomials satisfying the MIP's with respect to the matrices $\delta$ and $\eta$}\label{matrix:3}
        \begin{subfigure}[b]{0.3\textwidth}
        \includegraphics[width=\textwidth]{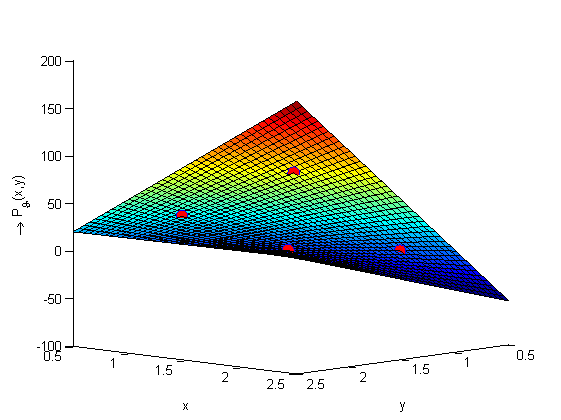}
        \caption{$P_{\vartheta} \in {\mathcal{P}_{2}^{2}}$}
        \label{fig:vartheta}
    \end{subfigure}
    ~ %add desired spacing between images, e. g. ~, \quad, \qquad, \hfill etc. 
    \begin{subfigure}[b]{0.3\textwidth}
        \includegraphics[width=\textwidth]{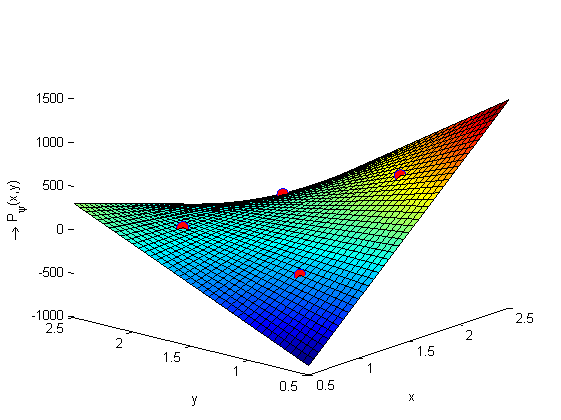}
        \caption{$P_{\psi} \in {\mathcal{P}_{2}^{2}}$}
        \label{fig:psi}
    \end{subfigure}
    \begin{subfigure}[b]{0.3\textwidth}
        \includegraphics[width=\textwidth]{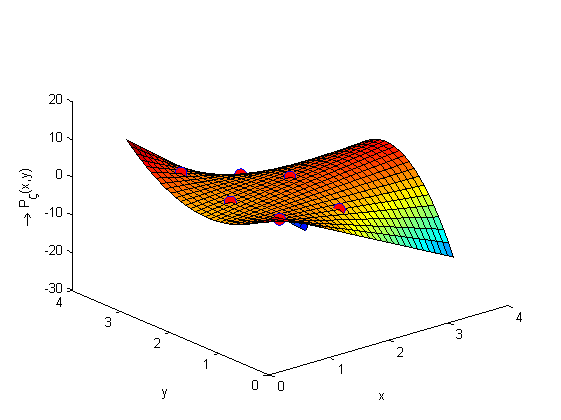}
        \caption{$P_{\zeta} \in {\mathcal{P}_{2}^{3}}$}
        \label{fig:zeta}
    \end{subfigure}
    ~ %add desired spacing between images, e. g. ~, \quad, \qquad, \hfill etc. 
    \begin{subfigure}[b]{0.3\textwidth}
        \includegraphics[width=\textwidth]{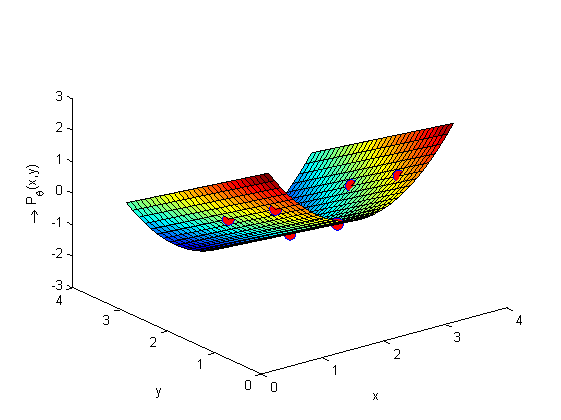}
        \caption{$P_{\theta} \in {\mathcal{P}_{3}^{2}}$}
        \label{fig:theta}
    \end{subfigure}
    \caption{Three-dimensional interpolating surfaces of the polynomials associated to the MIP's for the given matrices $\delta$, $\eta$, $\vartheta$, $\psi$, $\zeta$ and $\theta$ indicating the data points.}\label{matrix:6}
    \end{figure}
\end{example}

\begin{example}\label{ex:2}
Let $A_1 \in \mathbb{R}^{2 \times 3}$, $A_2 \in \mathbb{R}^{3 \times 2}$, $A_3 \in \mathbb{R}^{3 \times 3}$, and $A_4 \in \mathbb{R}^{3 \times 3}$ be four matrices given by 
$A_1 = \begin{pmatrix}		
		1 & -1 & 2\\
		-1 & 0 & -2
		    \end{pmatrix}$, $A_2 = \begin{pmatrix}	
		1 & -1\\
		-1 & 0\\
		2 & -2
		    \end{pmatrix}$, $A_3 = \begin{pmatrix}		
		1 & 2 & 3\\
		2 & 0 & 4\\
		3 & 4 & -1
		    \end{pmatrix}$, and $A_4 = \begin{pmatrix}		
		0 & 1 & 2\\
		-1 & 0 & 3\\
		-2 & -3 & 0
		    \end{pmatrix}$ respectively. Then, there exists unique $P_{A_1} \in \mathcal{P}_2^3$, $P_{A_2} \in \mathcal{P}_3^2$, $P_{A_3} \in \mathcal{P}_3^3$, and $P_{A_4} \in \mathcal{P}_3^3$ given as
\begin{eqnarray}
P_{A_1}(x,y) = -4xy^2+\frac{13}{2}y^2+15xy-13x-\frac{49}{2}y+21, \label{pol:A1}\\
%\end{flalign*}
%\begin{flalign*}
P_{A_2}(x,y) =  -4x^2y+\frac{13}{2}x^2+15xy-\frac{49}{2}x-13y+21, \label{pol:A2}
%\end{flalign*}
%\begin{flalign*}
\end{eqnarray}
\begin{eqnarray} \label{pol:A3}
P_{A_3}(x,y) =  -\frac{9}{2}x^2y^2+\frac{33}{2}x^2y+\frac{33}{2}xy^2-\frac{123}{2}xy-12x^2-12y^2\\+46x+46y-34, \nonumber
\end{eqnarray}
%\end{flalign*}
%\begin{flalign*}
\begin{eqnarray} \label{pol:A4}
\textrm{and } P_{A_4}(x,y) =  -x^2y+xy^2+x^2-y^2-2x+2y,
\end{eqnarray}
respectively. Clearly, $P_{A_1}(x,y) = P_{A_2}^{T}(x,y)$, $P_{A_3}(x,y) = P_{A_3}^{T}(x,y)$, and $P_{A_4}(x,y) = - P_{A_4}^{T}(x,y)$. Since $A_1=A_2^T$, $A_3=A_3^T$, and $A_4=-A_4^T$, where $A^T$ represents the transpose of the matrix $A$, i.e., the properties of transpose, symmetric and skew-symmetric matrices holds in the associated polynomial subspaces of two variables.

\begin{figure}[H]
    \centering
    \begin{subfigure}[b]{0.3\textwidth}
        \includegraphics[width=\textwidth]{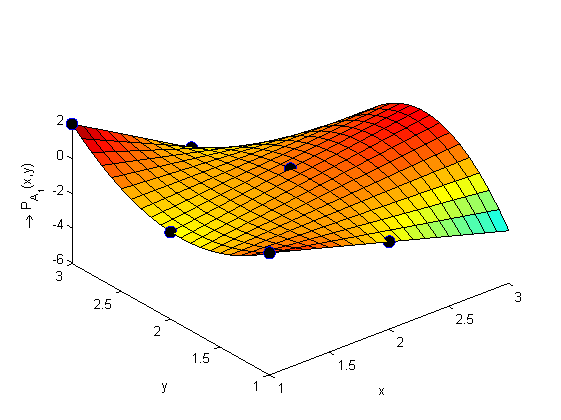}
        \caption{$P_{A_1} \in {\mathcal{P}_{2}^{3}}$}
        \label{fig:aone}
    \end{subfigure}
    ~ %add desired spacing between images, e. g. ~, \quad, \qquad, \hfill etc. 
    \begin{subfigure}[b]{0.3\textwidth}
        \includegraphics[width=\textwidth]{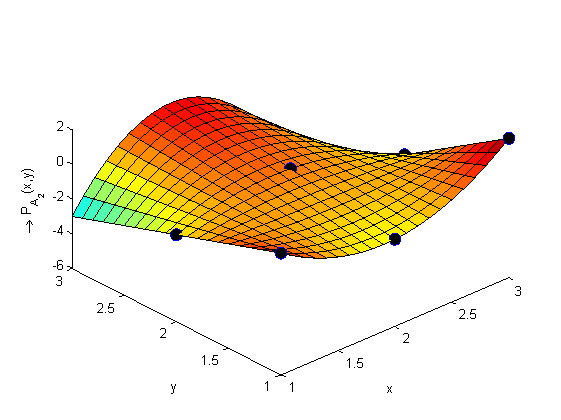}
        \caption{$P_{A_2} \in {\mathcal{P}_{3}^{2}}$}
        \label{fig:atwo}
    \end{subfigure}
    %\caption{Graphs of the polynomials satisfying the MIP's with respect to the matrices $\delta$ and $\eta$}\label{matrix:3}
        \begin{subfigure}[b]{0.3\textwidth}
        \includegraphics[width=\textwidth]{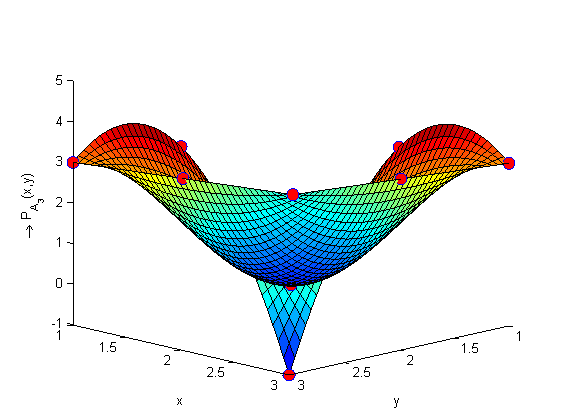}
        \caption{$P_{A_3} \in {\mathcal{P}_{3}^{3}}$}
        \label{fig:athree}
    \end{subfigure}
    ~ %add desired spacing between images, e. g. ~, \quad, \qquad, \hfill etc. 
    \begin{subfigure}[b]{0.3\textwidth}
        \includegraphics[width=\textwidth]{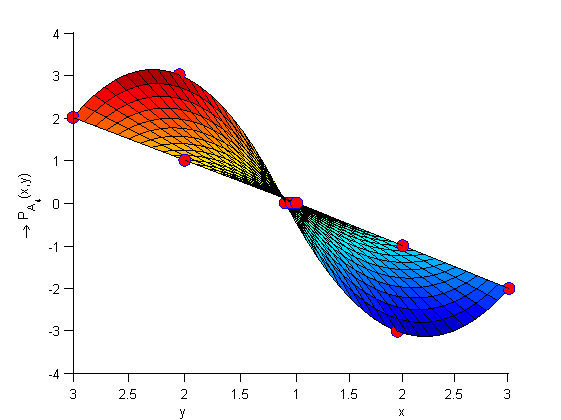}
        \caption{$P_{A_4} \in {\mathcal{P}_{3}^{3}}$}
        \label{fig:afour}
    \end{subfigure}
    ~ %add desired spacing between images, e. g. ~, \quad, \qquad, \hfill etc. 
    \begin{subfigure}[b]{0.3\textwidth}
        \includegraphics[width=\textwidth]{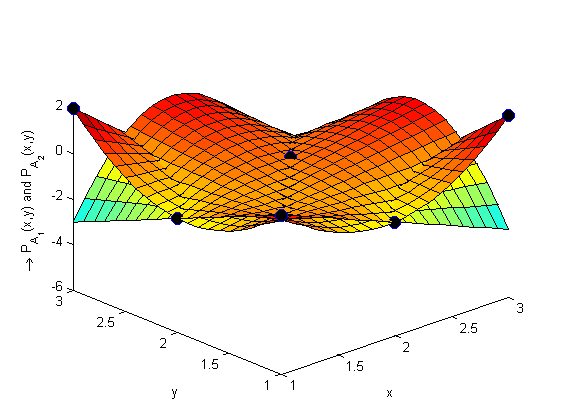}
        \caption{$P_{A_1} \in \mathcal{P}_{2}^{3}$ and $P_{A_2} \in \mathcal{P}_{3}^{2}$}
        \label{fig:aoneatwo}
    \end{subfigure}
    \caption{Three-dimensional interpolating surfaces validating geometric properties of the associated matrices $A_1$ $A_2$, $A_3$ and $A_4$ respectively.}\label{}
    \end{figure}
    
\end{example}

\begin{example}\label{ex:4}
Let $\tau = \begin{pmatrix}		
	    -1 & 2\\
		3 & -4
		\end{pmatrix}$
and $I_2 = \begin{pmatrix}		
	    1 & 0\\
		0 & 1
		    \end{pmatrix}$ be the identity matrix in $\mathbb{R}^{2 \times 2}$. Then, $\tau^2 = \begin{pmatrix}		
	    7 & -10\\
		-15 & 22
		\end{pmatrix}$ and there exists unique $P_\tau, P_{\tau^2}, \textrm{and } P_{I_2} \in \mathcal{P}_{2}^{2}$ given as $P_\tau(x,y) = -10xy+14x+13y-18$, $P_{\tau^2}(x,y) = 54xy-76x-71y+100$, and $P_{I_2}(x,y) = 2xy-3x-3y+5$ respectively. Here, $P_{\tau^2}(x,y) +5P_\tau(x,y) -2P_I(x,y)= 0$, i.e., 
\begin{equation}\label{charecteristric polynomial}
P_{\tau^2}(x,y) - tr(\tau).P_\tau(x,y) + det(\tau).P_{I_2}(x,y)= 0, 
\end{equation}
%        \begin{flalign*}
%	    \quad & P_\tau(x,y) = -10xy+14x+13y-18, \\ 
%		\quad & P_{\tau^2}(x,y) = 54xy-76x-71y+100, \\
%		\textrm{and } & P_I(x,y) = 2xy-3x-3y+5. &&
%		\end{flalign*}
i.e., the \textit{Cayley-Hamilton theorem} holds in the associated polynomial space $\mathcal{P}_{2}^{2}$. 

Additionally, since the polynomial $P_\tau(x,y) = -10xy+14x+13y-18$, is invertible (see example \ref{ex:invertible}). Let  $P_{\tau^{-1}} \in \mathcal{P}_{2}^{2}$ be the inverse polynomial. Then, using the equation \eqref{charecteristric polynomial} and definition \eqref{definition:dp product}, we get
\begin{equation}	P_{\tau^{-1}}(x,y)= \frac{1}{det(\tau)}[-P_{\tau}(x,y) + tr(\tau).P_I(x,y)]	=-\frac{1}{2}x-y+\frac{7}{2} \in \mathcal{P}_{2}^{2}.
	\end{equation}
Therefore, $\tau^{-1} = \begin{pmatrix}    
		2 & 1\\
		\frac{3}{2} & \frac{1}{2}
		\end{pmatrix}$ and it is the required inverse of the matrix $\tau \in \mathbb{R}^{2 \times 2}$.
\end{example}

\begin{example}\label{ex:5}
Let us consider the polynomial map $D_p:\mathbb{R}^{2 \times 2} \rightarrow \mathcal{P}_{2}^{2}$ defined as
\begin{equation}\label{linear:map2by2}
    \begin{aligned}
	   D_p{\left(\begin{bmatrix}		
	    a & b\\
		c & d
		\end{bmatrix} \right)} = (4a-2b & -2c+d) + (-2a+2b+c-d)x\\
		& -(2a-b-2c+d)y+(a-b-c+d)xy.
		\end{aligned}
\end{equation}
Since, ${D_p}{(kA+B)} = k{D_p}(A)+{D_p}(B)$ for all $A,B \in \mathbb{R}^{2 \times 2}$ and non-zero scalar $k$. The map \eqref{linear:map2by2} is linear. Again, $B_1=\left\{ \begin{bmatrix}		
	    1 & 0\\
		0 & 0
		\end{bmatrix}, \begin{bmatrix}		
	    0 & 1\\
		0 & 0
		\end{bmatrix}, \begin{bmatrix}		
	    0 & 0\\
		1 & 0
		\end{bmatrix}, \begin{bmatrix}		
	    0 & 0\\
		0 & 1
		\end{bmatrix}\right\}$ and $B_2={\{1,x,y,xy\}}$ are the standard bases of the spaces $\mathbb{R}^{2 \times 2}$ and $\mathcal{P}_{2}^{2}$ respectively, thus
\newline
\begin{flalign*}
	   D_p{\left(\begin{bmatrix}		
	    1 & 0\\
		0 & 0
		\end{bmatrix} \right)} = 4-2x-2y+xy,
		D_p{\left(\begin{bmatrix}		
	    0 & 1\\
		0 & 0
		\end{bmatrix} \right)} = -2+x+2y-xy, \\
D_p{\left(\begin{bmatrix}		
	    0 & 0\\
		1 & 0
		\end{bmatrix} \right)} =  -2+2x+y-xy \textrm{ and } D_p{\left(\begin{bmatrix}		
	    0 & 0\\
		0 & 1
		\end{bmatrix} \right)} = 1-x-y+xy.
\end{flalign*}
Therefore, the corresponding co-ordinate matrix with respect to the bases $B_1$ and $B_2$ is
\begin{equation}
	   [D_p]_{B_1}^{B_2} = \begin{bmatrix}		
	    4 & -2 & -2 & 1\\
		-2 & 1 & 2 & -1\\
		-2 & 2 & 1 & -1\\
		1 & -1 & -1 & 1
		\end{bmatrix}_{4 \times 4}.
\end{equation}
Here, $det\left([D_p]_{B_1}^{B_2} \right)=-1 \neq 0$, i.e., $ker(D_p)=\left\{0\right\}$.
\end{example}
\begin{note}
The inverse linear transformation $D_p^{-1}:\mathcal{P}_{2}^{2} \rightarrow \mathbb{R}^{2 \times 2}$ is given by
\begin{equation}
	   D_p^{-1}{\left(a+bx+cy+dxy \right)} = \begin{bmatrix}		
	    a+b+c+d & a+b+2c+2d\\
		a+2b+c+2d & a+2b+2c+4d
		\end{bmatrix}. \nonumber
\end{equation}
The corresponding co-ordinate matrix with respect to the bases $B_2$ and $B_1$ is
\begin{equation}
	   [D_p^{-1}]_{B_2}^{B_1} = \begin{bmatrix}		
	    1 & 1 & 1 & 1\\
		1 & 1 & 2 & 2\\
		1 & 2 & 1 & 2\\
		1 & 2 & 2 & 4
		\end{bmatrix}_{4 \times 4}.
\end{equation}
Clearly,
\begin{equation}
     [D_p]_{B_1}^{B_2}.[D_p^{-1}]_{B_2}^{B_1}=\begin{bmatrix}		
	    1 & 0 & 0 & 0\\
	    0 & 1 & 0 & 0\\
	    0 & 0 & 1 & 0\\
	    0 & 0 & 0 & 1
		\end{bmatrix}_{4 \times 4}=[D_p^{-1}]_{B_2}^{B_1}.[D_p]_{B_1}^{B_2}.  \nonumber
\end{equation}
\end{note}

\begin{example}\label{ex:6}
Let $A_1 \in \mathbb{R}^{2 \times 3}$, $A_2 \in \mathbb{R}^{3 \times 2}$, $A_3 \in \mathbb{R}^{3 \times 3}$, and $A_4 \in \mathbb{R}^{3 \times 3}$ be four matrices given in the example \ref{ex:2}. Then,  
${A_1}{A_2} = \begin{pmatrix}		
		6 & -5\\
		-5 & 5 
		    \end{pmatrix}$, 
${A_1}{A_3} = \begin{pmatrix}	
		5 & 10 & -3\\
		-7 & -10 & -1
		    \end{pmatrix}$, 
${A_1}{A_4} = \begin{pmatrix}	
		-3 & -5 & -1\\
		4 & 5 & -2
		    \end{pmatrix}$,
${A_2}{A_1} = \begin{pmatrix}		
		2 & -1 & 4\\
		-1 & 1 & -2\\
		4 & -2 & 8
		    \end{pmatrix}$, 
${A_3}{A_4} = \begin{pmatrix}		
		-8 & -8 & 8\\
		-8 & -10 & 4\\
		-2 & 6 & 18
		    \end{pmatrix}$,
and ${A_4}{A_1} = \begin{pmatrix}		
		 8 & 8 & 2\\
		8 & 10 & -6\\
		-8 & -4 & -18
		    \end{pmatrix}$.
Therefore, there exists unique $P_{{A_1}{A_2}} \in \mathcal{P}_2^2$, $P_{{A_1}{A_3}} \in \mathcal{P}_2^3$, $P_{{A_1}{A_4}} \in \mathcal{P}_2^3$, $P_{{A_2}{A_1}} \in \mathcal{P}_3^3$, $P_{{A_3}{A_4}} \in \mathcal{P}_3^3$, and $P_{{A_4}{A_3}} \in \mathcal{P}_3^3$ given as
\begin{flalign*}
P_{{A_1}{A_2}}(x,y) = & {} 21xy-32x-32y+49,\\
P_{{A_1}{A_3}}(x,y)  = & 15xy^2-24y^2-53xy+26x+85y-44,\\
P_{{A_1}{A_4}}(x,y) = & -7xy^2+10y^2+24xy-10x-35y+15,\\
P_{{A_2}{A_1}}(x,y) = & \frac{17}{2}x^2y^2-32x^2y-32xy^2+\frac{241}{2}xy+\frac{55}{2}x^2+\frac{55}{2}y^2-\frac{207}{2}x-\frac{207}{2}y+89,\\
P_{{A_3}{A_4}}(x,y) = & -3x^2y^2+15x^2y+9xy^2-47xy-9x^2+2y^2+29x+8y-12,\\
\textrm{and } P_{{A_4}{A_3}}(x,y) = & 3x^2y^2-9x^2y-15xy^2+47xy-2x^2+9y^2-8x-29y+12,&&
\end{flalign*}
respectively. The equations \eqref{pol:A1}, \eqref{pol:A2}, \eqref{pol:A3}, and \eqref{pol:A4} and the definition \eqref{definition:dp product} verify that $P_{{A_1}{A_2}}(x,y) = {P_{A_1}(x,y)} \otimes {P_{A_2}(x,y)}$, 
$P_{{A_1}{A_3}}(x,y)  = {P_{A_1}(x,y)} \otimes {P_{A_3}(x,y)}$,
$P_{{A_1}{A_4}}(x,y) = {P_{A_1}(x,y)} \otimes P_{A_4}(x,y)$,
$P_{{A_2}{A_1}}(x,y) = {P_{A_2}(x,y)} \otimes P_{A_1}(x,y)$,
$P_{{A_3}{A_4}}(x,y) = {P_{A_3}(x,y)} \otimes P_{A_4}(x,y)$,
\textrm{and } $P_{{A_4}{A_3}}(x,y) = {P_{A_4}(x,y)} \otimes P_{A_3}(x,y)$. 
\end{example}

\section{Concluding remarks}
In this paper, the well known $mn$-dimensional two-variate space of tensor-product polynomials, $\mathcal{P}_{m}^{n}$, of degree at most $(m-1,n-1)$, over the field $\mathbb{R}$ is considered. A theory of two-variable polynomials is developed in $\mathcal{P}_{m}^{n}$, parallel to the basic algebraic operations and geometric properties of the space $\mathbb{R}^{m \times n}$, by establishing the basic algebra and some algebraic properties with respect to the usual addition, scalar multiplication, and a newly defined algebraic operation `$\otimes$'. The space $\mathcal{P}_{m}^{n}$ form a subspace of two variable polynomials, of degree at most $m+n-2$ and is isomorphic to $\mathbb{R}^{m \times n}$ with respect to vector algebra structure under usual addition and scalar multiplication operation. Using counterexamples, it is verified that space $\mathcal{P}_{n}^{n}$ is non-commutative and possess zero-divisor property under the operation $\otimes$. It is proved that $\mathcal{P}_{m}^{n}$ is closed under the operation $\otimes$, if $m=n$. Some more algebraic properties of the polynomials in $\mathcal{P}_{m}^{n}$ are discussed under the operation $\otimes$ and it is proved that the associative and distributive law hold. Some definitions are included which define identity, inverse, eigenvalue, and power of the polynomials in $\mathcal{P}_{m}^{n}$ under the operation $\otimes$. In addition, the transpose, symmetry, skew-symmetry, and orthogonality of the polynomials in $\mathcal{P}_{n}^{n}$ is also defined. 

In the next section, the existence of the space $\mathcal{P}_{m}^{n}$ is established using tensor-product and univariate polynomial interpolation approach with respect to the MIP\eqref{mip} for a given matrix $(a_{ij}) \in \mathbb{R}^{m \times n}$. The poisedness of the MIP\eqref{mip} in $\mathcal{P}_{m}^{n}$ is also ensured. Further, three formulae are derived to construct the polynomial in $\mathcal{P}_{m}^{n}$ which satisfy the MIP\eqref{mip} for a given matrix $(a_{ij}) \in \mathbb{R}^{m \times n}$ using well-known univariate interpolation formulae. The method of the non-homogeneous system of linear equations is also remarked to obtain the coefficients of the required polynomial.

The next section focuses on some properties of the established polynomial map from $\mathbb{R}^{m \times n}$ to $\mathcal{P}_{m}^{n}$ on behalf of the construction formulae. It is achieved that the obtained polynomial map is invertible and linear, i.e., is an isomorphism. The respective invertible linear map is also defined. Finally, the multiplicative structure of the polynomials with respect to given matrices is discussed and it is proved that $\langle \mathbb{R}^{n \times n}, .\rangle$ and $\langle \mathcal{P}_{n}^{n}, \otimes  \rangle$ are isomorphic, i.e., the spaces $\mathbb{R}^{n \times n}$ and $\mathcal{P}_{n}^{n}$ are isomorphic with respect to the algebra structure. The existence of the identity polynomial in $\mathcal{P}_{m}^{n}$ under the operation $\otimes$ is also ensured which endorses that $\langle \mathcal{P}_{n}^{n},+,\otimes  \rangle$ form a ring with unity over $\mathbb{R}$.

In the last section, five examples are included to validate the results and for geometrical observations. The first example is validating the theorems \ref{thm:dp}, \ref{thm:construction by lagrange}, \ref{thm:construction by newton forwatd}, and \ref{thm:construction by newton backward} geometrically by surface diagrams in the figures \ref{matrix:6}. The second example shows that the polynomials in the associated subspaces preserve the geometric properties of the respective matrices such as transpose, symmetry, and skew-symmetry. In the third example, it is verified that the polynomials in the associated subspace $\mathcal{P}_{2}^{2}$ satisfy the result of Cayley-Hamilton theorem with respect to the given matrix in $\mathbb{R}^{2 \times 2}$ and the inverse of the given matrix is also obtained using polynomial algebra. In fourth example, the polynomial map from $\mathbb{R}^{2 \times 2}$ to $\mathcal{P}_{2}^{2}$ is defined using the result of the theorems \ref{thm:construction by lagrange}, \ref{thm:construction by newton forwatd} or \ref{thm:construction by newton backward}. It is verified that the defined map is linear as well as invertible. The inverse linear map is also included for the cross verification. Finally, last example verify that $P_{{A_i}{A_j}}(x,y)= P_{A_i}(x,y)\otimes P_{A_j}(x,y)$ (wherever defined) for $i,j=1,2,3,4$ with respect to the given matrices $A_1 \in \mathbb{R}^{2 \times 3}$, $A_2 \in \mathbb{R}^{3 \times 2}$, $A_3 \in \mathbb{R}^{3 \times 3}$, and $A_4 \in \mathbb{R}^{3 \times 3}$. Theoretical findings and numerical observations assert that space $\mathcal{P}_{m}^{n}$ preserves almost every algebraic and geometrical structure of the space $\mathbb{R}^{m \times n}$. 

The transformation of the finite order real matrices $a_{ij} \in \mathbb{R}^{m \times n}$ into $\mathcal{P}_{m}^{n}$, algebraic operations, respective algebra, algebraic and geometric structures, and isomorphism are building some additional tools for algebra, matrix theory, cryptography, and numerical linear algebra community. Considering the rapid growth in the development of the efficient algorithms for the polynomial computation (over matrices), the outcomes can also be beneficial for scientific and numerical computations. In addition, a matrix is a most common data structure to store the complete information of the objects, independent of the content of the data, in the various models of computer graphics, signal processing, image processing, computer-aided designs in control systems, thermal engineering to name a few. Therefore, the generation of the respective three-dimensional smooth interpolating surface (preserving the geometric nature) in the established subspace of two variable polynomials could play some important roles in the stages of analysis, interpretation, and applications.

%\textbf{SECTION-4 ADD HERE} 
%
%The invertible linear map given by \eqref{linearity}, satisfies the following properties:
%\begin{enumerate}[label=(\roman*)]
%\item $D_p(A+B) = D_{p}(A) + D_{p}(B)$ for all $A,B \in \mathbb{R}^{m \times n}$.
%\item $D_p(\alpha A) = \alpha D_{p}(A)$ for all $A \in \mathbb{R}^{m \times n}$, where $\alpha$ is a scalar.
%\item $D_p(A.B) = D_{p}(A)\otimes D_{p}(B)$ for all $A,B \in M_{n,n}$.
%\item $D_p(I_n) = I_n(x,y)$, where $I_n$ and $I_n(x,y)$ are the identity of $M_{n,n}$ and $\mathcal{P}_{n}^{n}$ respectively.
%\item $D_p(A^k) = ({D_p(A)})^k$  for all $A \in M_{n,n}$ and $k \in \mathbb{N}$.
%\item $D_p(A^{-1}) = ({D_p(A)})^{-1}$ for all invertible matrix $A \in M_{n,n}$.
%\end{enumerate}
%The corollary \eqref{coro:identity polynomial} insures that there exists identity polynomial in $\mathcal{P}_{n}^{n}$.

%The spaces $\mathcal{M}_{n,n}$ and $\mathcal{P}_{n}^{n}$ are isomorphic with respect to the algebra structure.

%\textbf{FUTURE WORK}
%\begin{conjecture}
%The MIP \eqref{mip} with respect to a given matrix $A \in \mathbb{R}^{m \times n}$ such that $mn=dim\Pi_{N}^{2}$, can not be poised in $\Pi_{N}^{2}$, for any $N \in \mathbb{N}$.
%\end{conjecture}

\appendix

\bibliography{sample}
\bibliographystyle{ieeetr}

\end{document}